\documentclass[twoside,11pt]{article}
\usepackage{jmlr2e}
\usepackage{array}
\usepackage{times}
\usepackage{amsmath}
\usepackage{url}
\usepackage{algorithm}
\usepackage{diagbox}
\usepackage{bbm}
\usepackage{bm}
\usepackage{graphicx}
\usepackage{subcaption}
\usepackage{epstopdf}
\usepackage{color}
\usepackage{algorithm}
\usepackage{algorithmic}
\usepackage{multirow}
\usepackage{rotating}
\usepackage{alphalph}

\newcommand{\bbb}[1]{\boldsymbol{\mathbf{#1}}}
\def\beq{\begin{eqnarray}}
\def\eeq{\end{eqnarray}}
\def\noi{\noindent}
\def\nn{\nonumber}
\def\la{\langle}
\def\ra{\rangle}

\newtheorem{assumption}{Assumption}
\def\ghs{\hspace{0.3cm}} 

\def\fourfigwid{0.4\textwidth}

\def\objimghei{0.15\textheight}

\firstpageno{1}
\let\cite\citep

\def\E{\mathbb{E}}
\def\P{\mathcal{P}}
\def\objimghei{0.17\textheight}
\def\objimgwid{1.0\textwidth}
\def\objimgheiB{0.17\textheight}
\begin{document}

\title{A Hybrid Method of Combinatorial Search and Coordinate Descent for Discrete Optimization}

\author{\name Ganzhao Yuan \email yuanganzhao@gmail.com \\
       \addr School of Data and Computer Science, Sun Yat-sen University (SYSU), P.R. China
       \AND
       \name Li Shen \email mathshenli@gmail.com \\
       \addr Tencent AI Lab, Shenzhen, P.R. China
       \AND
       \name Wei-Shi Zheng \email zhwshi@mail.sysu.edu.cn \\
       \addr School of Data and Computer Science, Sun Yat-sen University (SYSU), P.R. China}
\editor{}

\maketitle

\begin{abstract}

Discrete optimization is a central problem in mathematical optimization with a broad range of applications, among which binary optimization and sparse optimization are two common ones. However, these problems are NP-hard and thus difficult to solve in general. Combinatorial search methods such as branch-and-bound and exhaustive search find the global optimal solution but are confined to small-sized problems, while coordinate descent methods such as coordinate gradient descent are efficient but often suffer from poor local minima. In this paper, we consider a hybrid method that combines the effectiveness of combinatorial search and the efficiency of coordinate descent. Specifically, we consider random strategy or/and greedy strategy to select a subset of coordinates as the working set, and then perform global combinatorial search over the working set based on the original objective function. In addition, we provide some optimality analysis and convergence analysis for the proposed method. Our method finds stronger stationary points than existing methods. Finally, we demonstrate the efficacy of our method on some sparse optimization and binary optimization applications. As a result, our method achieves state-of-the-art performance in terms of accuracy. For example, our method generally outperforms the well-known orthogonal matching pursuit method in sparse optimization.

\end{abstract}

\section{Introduction}
In this paper, we mainly focus on the following nonconvex composite minimization problem (`$\triangleq$' means define):
\beq \label{eq:main}
 \min_{\bbb{x}\in \mathbb{R}^n}~F(\bbb{x}) \triangleq f(\bbb{x})+ h(\bbb{x})
\eeq
\noi where $f(\bbb{x})$ is a smooth convex function with its gradient being $L$-Lipschitz continuous, and $h(\bbb{x})$ is a piecewise separable function with $h(\bbb{x})=\sum_{i}^n h_i(\bbb{x}_i)$. We consider two cases for $h(\cdot)$:
\beq
\textstyle&\{ h_{\text{binary}}(\bbb{x}) \triangleq I_{\Psi}(\bbb{x}),~\Psi \triangleq \{-1,1\}^n \} ~~~~~\nn\\
\text{or}~~&\{ h_{\text{sparse}}(\bbb{x}) \triangleq \lambda \|\bbb{x}\|_0 + I_{\Omega}(\bbb{x}),~\Omega \triangleq [-\rho \bbb{1},\rho \bbb{1} ]\},\nn
\eeq
\noi where $I_{\Psi}(\cdot)$ is an indicator function on $\Psi$ with $I_{\Psi}(\bbb{x})={\tiny \left\{
                                                                   \begin{array}{ll}
                                                                     0, & \hbox{$\bbb{x}\in\Psi$;} \\
                                                                     +\infty, & \hbox{$\bbb{x}\notin\Psi$.}
                                                                   \end{array}
                                                                 \right.}$, $\|\cdot\|_0$ is a function that counts the number of nonzero elements in a vector, $\lambda$ and $\rho$ are strictly positive constants. When $h \triangleq h_{\text{binary}}$, (\ref{eq:main}) refers to the binary optimization problem; when $h \triangleq h_{\text{sparse}}$, (\ref{eq:main}) corresponds to the sparse optimization problem.

Binary optimization and sparse optimization capture a variety of applications of interest in both machine learning and computer vision, including binary hashing \cite{WangLKC16,Wang2017}, dense subgraph discovery \cite{yuan2013truncated,YuanG17aaai}, Markov random fields \cite{Boykov2001}, compressive sensing \cite{candes2005decoding,Donoho06}, sparse coding \cite{aharon2006img,BaoJQS16}, subspace clustering \cite{elhamifar2013sparse}, to name only a few. In addition, binary optimization and sparse optimization are closely related to each other. A binary optimization problem can be reformulated as a sparse optimization problem using the fact that \cite{YuanG17aaai}: $\bbb{x} \in \{-1,1\}^n \Leftrightarrow \|\bbb{x}-\bbb{1}\|_0 + \|\bbb{x}+\bbb{1}\|_0\leq n$, and the reverse is also true using the variational reformulation of $\ell_0$ pseudo-norm \cite{bienstock1996computational}: $\forall \|\bbb{x}\|_{\infty} \leq \rho,~\|\bbb{x}\|_0 = \min ~\la \bbb{1},\bbb{v} \ra,~s.t.~\bbb{v} \in \{0,1\}^n,~|\bbb{x}|\leq \rho \bbb{v}$. There are generally four classes of methods for solving the binary or sparse optimization problem in the literature, which we present below.

\textbf{Relaxed Approximation Method}. One popular method to solve (\ref{eq:main}) is convex or nonconvex relaxed approximation method. Box constrained relaxation, semi-definite programming relaxation, and spherical relaxation are often used for solving binary optimization problems, while $\ell_1$ norm, top-$k$ norm, Schatten $\ell_p$ norm, re-weighted $\ell_1$ norm, capped $\ell_1$ norm, half quadratic function, and many others are often used for solving sparse optimization problems. It is generally believed that nonconvex methods often achieve better accuracy than the convex counterparts. Despite the merits, this class of methods fails to directly control the sparsity or binary property of the solution.

\textbf{Greedy Pursuit Method}. This method is often used to solve cardinality constrained discrete optimization problems. For sparse optimization, this method greedily selects at each step one atom of the variables which have some desirable benefits \cite{tropp2007signal,dai2009subspace,needell2010signal,blumensath2008gradient,needell2009cosamp}. It has a monotonically decreasing property and achieves optimality guarantees in some situations, but it is limited to solving problems with smooth objective functions (typically the square function). For binary optimization, this method is strongly related to submodular optimization as minimizing a set function can be reformulated as a binary optimization problem \cite{CalinescuCPV11}.


\textbf{Combinatorial Search Method}. Combinatorial search method \cite{conforti2014integer} is typically concerned with problems that are NP-hard. A naive method is exhaustive search (a.k.a generate and test method). It systematically enumerates all possible candidates for the solution and pick the best candidate corresponding to the lowest objective value. The cutting plane method solves the convex linear programming relaxation and adds linear constraints to drive the solution towards binary variables, while the branch-and-cut method performs branches and applies cuts at the nodes of the tree having a lower bound that is worse than the current solution. Although in some cases these two methods converge without much effort, in the worse case they end up solving all $2^n$ convex subproblems.


\bbb{Proximal Point Method}. Based on the current gradient $\nabla f(\bbb{x}^k)$, proximal point method \cite{beck2013sparsity,lu2014iterative,jain2014iterative,nguyen2014linear,patrascu2015efficient,PatrascuN15,li2016nonconvex} iteratively performs a gradient update followed by a proximal operation: $\bbb{x}^{k+1}= \text{prox}_{\gamma h}(\bbb{x}^k - \gamma \nabla f(\bbb{x}^k))$. Here the proximal operator $\text{prox}_{\tilde{h}}(\bbb{a}) = \arg \min_{\bbb{x}} ~\tfrac{1}{2}\|\bbb{x}-\bbb{a}\|_2^2 + \tilde{h}(\bbb{x})$ can be evaluated analytically, and $\gamma={1}/{L}$ is the step size with $L$ being the Lipschitz constant. This method is closely related to (block) coordinate descent \cite{nesterov2012efficiency, chang2008coordinate,breheny2011coordinate,de2016fast,razaviyayn2013unified,beck2013convergence,hong2013iteration,lu2015complexity,xu2013block} in the literature. Due to its simplicity, many strategies (e.g., variance reduction \cite{johnson2013accelerating,xiao2014proximal,chen2016accelerated}, asynchronous parallelism \cite{liu2015asynchronous,recht2011hogwild}, and non-uniform sampling \cite{zhang2016accelerated}) have been proposed to accelerate proximal point method. However, existing works use a scalar step size and solve a first-order majorization/surrogate function via closed form updates. Since problem (\ref{eq:main}) is nonconvex, such a simple majorization function may not necessarily be a good approximation for the original problem.


Compared to the four existing solutions mentioned above, our method has the following four merits. \bbb{(i)} It can directly control the sparsity or binary property of the solution. \bbb{(ii)} It is a greedy coordinate descent algorithm \footnote{This is in contract with greedy pursuit method where the solutions must be initialized to zero and may cause divergence when being incorporated to solve bilinear matrix factorization \cite{BaoJQS16}.}. \bbb{(iii)} It leverages the effectiveness of combinatorial search. \bbb{(iv)} It significantly outperforms proximal point method and inherits its computational advantages.




The contributions of this paper are three-fold. \bbb{(i)} Algorithmically, we introduce a novel hybrid method (denoted as HYBRID) for sparse or binary optimization which combines the effectiveness of combinatorial search and the efficiency of coordinate descent (See Section \ref{sect:proposed}). \bbb{(ii)} Theoretically, we establish the optimality hierarchy of our proposed algorithm and show that it always finds a stronger stationary point than existing methods (See Section \ref{sect:hierarchy}). In addition, we prove the global convergence and convergence rate of the proposed algorithm (See Section \ref{sect:convergence}). \bbb{(iii)} Empirically, we have conducted extensive experiments on some some binary optimization and sparse optimization tasks to show the superiority of our method (See Section \ref{sect:exp}).

\section{Proposed Algorithm} \label{sect:proposed}

This section presents our hybrid method for solving the optimization problem in (\ref{eq:main}). Our algorithm is an iterative procedure. In every iteration, the index set of variables is separated to two sets $B$ and $N$, where $B$ is the working set. We fix the variables corresponding to $N$, while minimize a sub-problem on variables corresponding to $B$. We use $\bbb{x}_B$ to denote the sub-vector of $\bbb{x}$ indexed by $B$. The proposed method is summarized in Algorithm \ref{algo:main}.

\begin{algorithm}[!t]
\caption{ {\bf A Hybrid Approach for Sparse or Binary Optimization} }
\begin{algorithmic}[1]
  \STATE Input: the size of the working set $k$, an initial feasible solution $\bbb{x}^0$. Set $t=0$.
\WHILE{not converge}
\STATE (S1) Employ some strategy to find a working set $B$ of size $k$. Denote $N \triangleq   \{1,...,n\}\setminus B$.
  \STATE (S2) Solve the following subproblem \emph{globally} using combinatorial search:
  \beq \label{eq:subprob}
  \begin{split}
\textstyle  \bbb{x}^{t+1} \Leftarrow \arg \min_{\bbb{z}}~f(\bbb{z}) + h(\bbb{z}) + \tfrac{\theta}{2} \|\bbb{z}-\bbb{x}^{t}\|^2,\\
~s.t.~\bbb{z}_{N} = \bbb{x}^t_{N}~~~~~~~~~~~~~~~~~~~~~~~~
  \end{split}
\eeq
\STATE (S3) Increment $t$ by 1
\ENDWHILE
\end{algorithmic}
\label{algo:main}
\end{algorithm}

At first glance, Algorithm \ref{algo:main} might seem to be merely a (block) coordinate descent algorithm \cite{tseng2009coordinate} applied to (\ref{eq:main}). However, it has some interesting properties that are worth commenting on.

\noi \bbb{$\bullet$ Two New Strategies.} \bbb{(i)} Instead of using majorization techniques for optimizing the block of variable, we consider minimizing the original objective function. Although the subproblem is NP-hard and admits no closed form solution, we can use an exhaustive search to solve it exactly. \bbb{(ii)} We consider a proximal point strategy for the subproblem. This is to guarantee sufficient descent condition of the objective function and global convergence of Algorithm \ref{algo:main}(refer to Theorem \ref{theorem:convergence}).


\noi \bbb{$\bullet$ Solving the Subproblem Globally.} The subproblem in (\ref{eq:subprob}) essentially contains $k$ unknown decision variables and can be solved exactly within sub-exponential time $\mathcal{O}(2^k)$. For both binary optimization and sparse optimization, problem (\ref{eq:subprob}) can be reformulated as a integer/mixed-integer optimization problem and solved by global optimization solvers such as CPLEX or Gurobi. For simplicity, we consider a simple exhaustive search to solve it. Specifically, for every coordinate of the $k$-dimensional subproblem, it has two states, i.e., zero/nonzero. We systematically enumerate the full binary tree to obtain all possible candidate solutions and then pick the best one that leads to the lowest objective value as the optimal solution \footnote{We take $f(\bbb{x}) \triangleq \tfrac{1}{2}\bbb{x}^T\bbb{Qx} + \la \bbb{x},\bbb{p} \ra$ and $h(\bbb{x}) \triangleq h_{\text{sparse}}(\bbb{x})$ for example, where $\bbb{Q}\in\mathbb{R}^{n\times n}$ and $\bbb{p}\in \mathbb{R}^n$ are given. Problem (\ref{eq:subprob}) is equivalent to the following small-sized optimization problem: $\bbb{v}^*=\arg\min_{\bbb{v} \in \mathbb{R}^{k}}~\phi(\bbb{\bbb{v}}) \triangleq \frac{1}{2}\bbb{v}^T \bbb{Q}_{B,B}\bbb{v} + \la \bbb{v}, \bbb{Q}_{B,N}\bbb{x}^t_N + \bbb{p}_{B} \ra + h(\bbb{v}) + \tfrac{\theta}{2} \|\bbb{v}-\bbb{x}_B^t\|_2^2$ with $\bbb{x}^{t+1}_B=\bbb{v}^*$ and $\bbb{x}^{t+1}_N=\bbb{x}_N^t$. We consider minimizing the following problem $\min_{\bbb{v} \in \mathbb{R}^{k}}~\phi(\bbb{\bbb{v}}),~s.t.~\bbb{v}_K=\bbb{0}$, where $K$ contains $\sum_{i=0}^k C_k^i$ possible choices for the coordinates.  }.  




\noi \bbb{$\bullet$ Finding a Working Set.} We observe that it contains $C_n^k$ possible combinations of choice for the working set. One may use a cyclic strategy to alternatingly select all the choices of the working set. However, past results show that coordinate gradient method results in faster convergence when the working set in selected in an arbitrary order \cite{hsieh2008dual} or in a greedy manner \cite{tseng2009coordinate,hsieh2011fast}. This inspires us to use random strategy and greedy strategy for finding the working set in Algorithm \ref{algo:main}. We remark that the combination of the two strategies is preferred in practice.

\boxed{\text{Random strategy}.} We uniformly select one combination (which contains $k$ coordinates) from the whole working set of size $C_n^k$. One remarkable benefit of this strategy is that our algorithm is ensured to find the block-$k$ stationary point in expectation.

\boxed{\text{Greedy strategy}.} We use the following variable selection strategy for sparse optimization, and similar strategy can be directly applied to binary optimization. Generally speaking, we pick top-$k$ coordinates that lead to the greatest descent when one variable is changed and the rest variables are fixed based on the current solution $\bbb{x}^t$. We denote $I\triangleq \{i: \bbb{x}^t_i=0\}$ and $J\triangleq \{j: \bbb{x}^t_j \neq 0\}$. We expect the working set is balanced and pick $k/2$ coordinates from $I$ and $k/2$ coordinates from $J$. For $I$, we solve a one-variable subproblem to compute the possible decrease for all $i \in I$ of $\bbb{x}^t$ when changing from zero to nonzero:
\beq
 \forall i =1,...,|I|,~ \bbb{c}_{i} = \min_{\alpha} F(\bbb{x}^t + \alpha \bbb{e}_i) - F(\bbb{x}^t). \nn
  \eeq
  \noi For $J$, we compute the decrease for each coordinate $j \in J$ of $\bbb{x}^t$ when changing from nonzero to exactly zero:
  \beq
 \forall j =1,...,|J|,~\bbb{d}_{j} = F(\bbb{x}^t + \alpha \bbb{e}_j) - F(\bbb{x}^t),~\alpha = \bbb{x}^t_j. \nn
  \eeq
  \noi Here $\bbb{e}_i$ is a unit vector with a $1$ in the $i$th entry and $0$ in all other entries. Assuming that $k$ is an even number, we sort the vectors $\bbb{c}$ and $\bbb{d}$ in increasing order and then pick top-$(k/2)$ coordinates form $I$ and top-$(k/2)$ coordinates from $J$ as the working set \footnote{Continuing our previous example with $f(\bbb{x}) \triangleq \tfrac{1}{2}\bbb{x}^T\bbb{Qx} + \la \bbb{x},\bbb{p} \ra$ and $h(\bbb{x}) \triangleq h_{\text{sparse}}(\bbb{x})$, the vectors $\bbb{c} \in \mathbb{R}^{|I|}$ and $\bbb{d}\in \mathbb{R}^{|J|}$ can be further computed as $\bbb{c}_i = \alpha (\bbb{Qx}+\bbb{p})_i + 0.5 \alpha ^2 \bbb{Q}_{i,i} + \lambda$ and $\bbb{d}_j = \alpha (\bbb{Qx}+\bbb{p})_j + 0.5 \alpha ^2 \bbb{Q}_{j,j} - \lambda$, respectively.}. If either $|I|<k/2$ or $|J|<k/2$, one can pick the whole set of $I$ or $J$ as a part of the working set.


\noi \bbb{$\bullet$ Extensions to Cardinality Constrained Problems.} In many applications, it is desirable to directly control the cardinality of the solution using the following constraints: 
\beq
\{ h_{\text{binary-c}}(\bbb{x}) \triangleq I_{\Upsilon}(\bbb{x}),~\Upsilon \triangleq \{\bbb{x}~|~\bbb{x}\in \{0,1\}^n,~\bbb{x}^T\bbb{1}=s\} \}~\nn\\
\text{or}~\{ h_{\text{sparse-c}}(\bbb{x}) \triangleq I_{\Phi}(\bbb{x}),~\Phi \triangleq \{\bbb{x}~|~\|\bbb{x}\|_0\leq s \}\}. ~~~~~~~~~~\nn
\eeq
\noi The proposed block coordinate method (including the exhaustive search algorithm and the working set selection strategies) can still be applied even when $h(\bbb{x})$ contains one non-separable constraint. What one needs is to ensure that the solution $\bbb{x}^{t}$ is a feasible solution for all $t=0,1,...\infty$. This is similar to the prior work of \cite{necoara2013random}.

\section{Optimality Analysis}\label{sect:hierarchy}

This section provides some optimality analysis for our method. In the sequel, we present some necessary optimal conditions for (\ref{eq:main}). Since the block-$k$ optimality is novel in this paper, it is necessary to clarify its relations with existing optimality conditions formally. We use $\breve{\bbb{x}},~\grave{\bbb{x}}$, and $\bar{\bbb{x}}$ to denote a basic stationary point, an $L$-stationary point, and a block-$k$ stationary point, respectively.

\begin{definition} \label{def:basic}
(Basic Stationary Point) A solution $\breve{\bbb{x}}$ is called a basic stationary point if the following holds. $h \triangleq h_{\text{binary}}: \breve{\bbb{x}} \in \{-1,+1\}^n;~h \triangleq h_{\text{sparse}}:\breve{\bbb{x}}_S = \arg\min_{ \bbb{z} \in[-\rho\bbb{1},\rho\bbb{1}]}\tfrac{L}{2}\|\bbb{z} - ( \breve{\bbb{x}} - \nabla f(\breve{\bbb{x}})/L)_S \|_2^2$, where $S\triangleq \{i| \breve{\bbb{x}}_i \neq 0\}$.
\end{definition}

\noi \textbf{Remarks:} For binary optimization, any binary solution is a basic stationary point. For sparse optimization, basic stationary point states that the solution achieves its global optimality when the support set is restricted. One remarkable feature of the basic stationary condition is that the solution set is enumerable and its size is $2^n$. It makes it possible to validate whether a solution is optimal for the original discrete optimization problem.

\begin{definition} \label{def:L:stationary}
($L$-Stationary Point) A solution $\grave{\bbb{x}}$ is an $L$-stationary point if it holds that: $\grave{\bbb{x}} = \arg\min_{\bbb{x}}~g(\bbb{x},\grave{\bbb{x}})+ h(\bbb{x})$ with $g(\bbb{x},\bbb{z})\triangleq f(\bbb{z}) + \la \nabla f(\bbb{z}),\bbb{x} - \bbb{z} \ra + \frac{L}{2}\|\bbb{x} - \bbb{z}\|_2^2,~\forall \bbb{x},~\bbb{z},~f(\bbb{x})\leq  g(\bbb{x},\bbb{z})$.
\end{definition}

\noi \textbf{Remarks:} This is the well-known proximal thresholding operator \cite{beck2013sparsity}. Although it has a closed-form solution, this simple surrogate function may not be a good majorization/surrogate function for the non-convex problem.

\begin{definition} \label{def:block:k}
 (Block-$k$ Stationary Point) A solution $\bar{\bbb{x}}$ is a block-$k$ stationary point if it holds that:
\beq
\bar{\bbb{x}} \in \arg\min_{\bbb{z}\in\mathbb{R}^n}~\P(\bbb{z};\bar{\bbb{x}},B)\triangleq \{F(\bbb{z}),~s.t.~\bbb{z}_N=\bar{\bbb{x}}_N\},
\eeq
\noi for all $|B|=k,~N \triangleq \{1,...,n\}\setminus B$.
\end{definition}

\noi\textbf{Remarks:} Block-$k$ stationary point is novel in this paper. One remarkable feature of this concept is that it involves solving a small-sized NP-hard problem which can be tackled by some practical global optimization method.

The following proposition states the relations between the three types of stationary point.
\begin{proposition} \label{proposition:hierarchy}
\textbf{Proof of the Hierarchy between the Necessary Optimality Conditions.} We have the following optimality hierarchy: $\boxed{\text{Basic Stat. Point}} \overset{(1)}{\Leftarrow} \boxed{ \text{$L$-Stat. Point}} \overset{(2)}{\Leftarrow}  \boxed{  \text{Block-$1$ Stat. Point}} \overset{(3)}{\Leftarrow}  \boxed{  \text{Block-$2$ Stat. Point}} \overset{}{\Leftarrow}  ... \overset{}{\Leftarrow}  \boxed{  \text{Block-$n$ Stat. Point}}   \overset{(4)}{\Leftrightarrow}  \boxed{ \text{Optimal Point}}$.

\begin{proof}
For notational convenience, we denote $\Pi(\bbb{a}) \triangleq \min(\rho\bbb{1}, \max(-\rho\bbb{1},\bbb{a}))$.

(1) We now prove that an $L$-stationary point $\grave{\bbb{x}}$ is also a basic stationary point $\breve{\bbb{x}}$. When $h \triangleq h_{\text{binary}}$, this conclusion clearly holds. We now consider $h \triangleq h_{\text{sparse}}$. Observing that the problem (see Definition \ref{def:basic}) is separable, we have the following closed-form solution: $\breve{\bbb{x}}_S = \Pi(\breve{\bbb{x}}_S - (\nabla f(\breve{\bbb{x}}))_S/L)$. For the problem in Definition \ref{def:L:stationary}, we have the following closed form solution for $\grave{\bbb{x}}$:
{\small \beq
\grave{\bbb{x}}_i = \left\{
              \begin{array}{ll}
                \Pi(\grave{\bbb{x}}_i - \nabla_i f(\grave{\bbb{x}})/L), & \hbox{$(\grave{\bbb{x}}_i - \nabla_i f(\grave{\bbb{x}})/L)^2 > 2 \lambda/L$;} \\
                0, & \hbox{$\text{else}.$}
              \end{array}
            \right.\nn
 \eeq }{\normalsize}\noi Clearly, the latter formulation implies the former one. Moreover, it is not hard to notice that for all $i$ that $\grave{\bbb{x}}_i=0$, we have $|\nabla_i f(\grave{\bbb{x}})| \leq \sqrt{2 \lambda L}$, and for all $j$ that $\grave{\bbb{x}}_j \neq 0$, we have $\nabla_j f(\grave{\bbb{x}}) = 0$ and $ |\grave{\bbb{x}}_j | \geq \min(\rho,\sqrt{2 \lambda/L})$.

(2) Note that the convex objective function $f(\cdot)$ has coordinate-wise Lipschitz continuous gradient with constant $\bbb{s}_i$. For all $\bbb{x}\in\mathbb{R}^n,~\delta\in\mathbb{R},~i=1,2,...n$, it holds that: \cite{nesterov2012efficiency}
$$f(\bbb{x}+ \delta \bbb{e}_i) \leq Q_i(\bbb{x},\delta)\triangleq f(\bbb{x}) + \la \nabla_i f(\bbb{x}),\delta\bbb{e}_i \ra + \frac{\bbb{s}_i}{2}\| \delta\bbb{e}_i\|_2^2,$$
\noi where $\bbb{e}_i \in \mathbb{R}^n$ a unit vector with a 1 in the $i$th entry and 0 in all other entries. Any block-1 stationary point must satisfy the following relation:
$$0 \in \arg \min_{\delta}~Q_i(\bar{\bbb{x}},\delta)+ h_i(\bar{\bbb{x}}_i + \delta )$$
\noi for all $i$. For $h=h_{\text{binary}}$, we have the following result:
\beq
\bar{\bbb{x}}_i = {\left\{
              \begin{array}{ll}
                1, & \hbox{$\bar{\bbb{x}}_i - \nabla_i f(\bar{\bbb{x}})/\bbb{s}_i > 0$;} \\
                -1, & \hbox{$\text{else}$.}
              \end{array}
            \right\}}.\nn
\eeq
\noi For $h=h_{\text{sparse}}$, we have the following result:
{\small\beq
\bar{\bbb{x}}_i = {  \left\{
              \begin{array}{ll}
                \Pi(\bar{\bbb{x}}_i - \nabla_i f(\bar{\bbb{x}})/\bbb{s}_i  ), & \hbox{$(\bar{\bbb{x}}_i - \nabla_i f(\bar{\bbb{x}})/\bbb{s}_{i})^2 > \frac{2\lambda}{\bbb{s}_i}$;} \\
                0, & \hbox{$\text{else}$.}
              \end{array}
            \right\}}.\nn
\eeq}{\normalsize}\noi Since $\bbb{s}_i \leq L$ for all $i$, we conclude that block-1 stationary point implies $L$-stationary point.

(3) We now show that block-$k_1$ stationary point implies block-$k_2$ stationary point when $k_1\geq k_2$. Note that to guarantee block-$k$ stationary condition, one need to solve the problem in Definition \ref{def:block:k} for $\sum_{i=0}^k C_{n}^k$ times, i.e. all the combination which is at most of size $k$. Clearly, when $k_1\geq k_2$, the subproblems for block-$k_2$ stationary point is a subset of those of block-$k_1$ stationary point.

(4) It is clearly that any block-$n$ stationary point is also the global optimal soluiton.

\end{proof}

\end{proposition}

\noi\textbf{Remarks:} It is worthwhile to point out that the seminal work of \cite{beck2013sparsity} also presents an optimality condition for sparse optimization. However, our block-$k$ condition is stronger than their coordinate-wise optimality since their optimal condition corresponds to $k=1$ in our optimality condition framework.

\textbf{A Running Example.} We consider the quadratic optimization problem $\min_{\bbb{x}\in\mathbb{R}^n}~\tfrac{1}{2}\bbb{x}^T\bbb{Qx} + \bbb{x}^T\bbb{p} + h(\bbb{x})$ where $n=6$,~$\bbb{Q}=\bbb{c}\bbb{c}^T+\bbb{I}$,~$\bbb{p}=\bbb{1}$,~$\bbb{c}=[1~2~3~4~5~6]^T$ when $h(\cdot)\triangleq h_{\text{sparse}}(\bbb{x})$ and $h(\cdot)\triangleq h_{\text{binary}}(\bbb{x})$. The parameters for $h_{\text{sparse}}(\bbb{x})$ are set to $\lambda=0.01$ and $\rho = +\infty$. The stationary point distribution on this example can be found in Table \ref{tab:optimality}. This problem contains $\sum_{i=0}^6 C_6^i = 64$ stationary points. There are 56 and 58 local minimizers satisfying the $L$-stationary condition for binary optimization and sparse optimization problem, respectively. There are 9 and 11 local minimizers satisfying the block-1 stationary condition. Moreover, as $k$ becomes large, the newly introduced type of local minimizers (i.e. block-$k$ stationary point) become more restricted in the sense that they have small number of stationary points.

\begin{table}[!h]
\fontsize{7.2}{7.2}\selectfont
\centering
\scalebox{1}{\begin{tabular}{|p{1.2cm}|p{1.2cm}|p{1.2cm}|p{1.2cm}|p{1.2cm}|p{1.2cm}|p{1.2cm}|p{1.2cm}|p{1.3cm}|}
  \hline
 & {\footnotesize  \text{Basic}-\text{Stat.}} & {\footnotesize L-\text{Stat.}} & {\footnotesize  \text{Block}-1 \text{Stat.}} & {\footnotesize \text{Block}-2 \text{Stat.} }&  {\footnotesize \text{Block}-3 \text{Stat.}}& {\footnotesize  \text{Block}-4 \text{Stat.}}&  {\footnotesize \text{Block}-5 \text{Stat.}}&  {\footnotesize \text{Block}-6 \text{Stat.}}\\
  \hline
$h \triangleq h_{\text{binary}}$ &  64 & 56 &  9 & 3 &  1 &  1&  1&  1 \\
  \hline
$h \triangleq h_{\text{sparse}}$ &  64 & 58 &  11 & 2 &  1 &  1&  1&  1 \\
  \hline
\end{tabular}}
\caption{Number of points satisfying optimality conditions.}\label{tab:optimality}
\end{table}

\section{Convergence Analysis} \label{sect:convergence}

This section provides some convergence analysis for Algorithm \ref{algo:main}. We assume that the working set of size $k$ is selected uniformly.

We now prove the global convergence for general discrete optimization. We provide the convergence rate of our algorithm in the \textbf{appendix}.

%





\begin{theorem} \label{theorem:convergence}

\textbf{Proof of Global Convergence Properties.} Letting $\bbb{x}^t$ be the sequence generated by Algorithm \ref{algo:main}, we have the following results.

\bbb{(i)} For $h \triangleq h_{\text{binary}}$, $h \triangleq h_{\text{sparse}}$, $h \triangleq h_{\text{binary-c}}$, and $h \triangleq h_{\text{sparse-c}}$, it holds that:
\beq\label{eq:suff:dec}
&\textstyle F(\bbb{x}^{t+1}) - F(\bbb{x}^t) \leq - \tfrac{\theta}{2}\|\bbb{x}^{t+1}-\bbb{x}^{t}\|^2,\nn\\
&\textstyle\lim_{t\rightarrow \infty} \E[\|\bbb{x}^{t+1} - \bbb{x}^t\|] = 0.\nn
\eeq

\bbb{(ii)} As $t\rightarrow \infty$, $\bbb{x}^t$ converges to the block-$k$ stationary point $\bar{\bbb{x}}$ of (\ref{eq:main}) in expectation.

\bbb{(iii)} For $h\triangleq h_{\text{binary}}$, it holds that $\E[\|\bbb{x}^{t+1}-\bbb{x}^t\|_2^2] \geq \frac{2k}{n}$. The solution changes at most $\frac{n [F(\bbb{x}^0) - F(\bar{\bbb{x}})]}{2k \theta}$ times in expectation for finding a block-$k$ stationary point.

\bbb{(iv)} For $h\triangleq h_{\text{sparse}}$, it holds that $|\bbb{x}_i^{t}|\geq \delta$ for all $i$ with $\bbb{\bbb{x}}_i^t\neq 0$, where $\delta\triangleq \min(\rho,\sqrt{ \mathstrut 2\lambda/(\theta + L)},\min(|\bbb{x}^0|))$. In addition, whenever $\bbb{x}^{t+1} \neq \bbb{x}^t$, we have: $\E[\|\bbb{x}^{t+1}-\bbb{x}^t\|_2^2] \geq \frac{k\delta^2}{n}$, the objective value is decreased at least by $D$, and the solution changes at most ${\bar{J}}$ times in expectation for finding a block-$k$ stationary point $\bar{\bbb{x}}$. Here $D$ and ${\bar{J}}$ are defined as
\beq \label{eq:DD}
D \triangleq \frac{k \theta\delta^2}{2n},~~{\bar{J}} \triangleq \frac{F(\bbb{x}^0) - F(\bar{\bbb{x}})}{D}
\eeq

\begin{proof}

\bbb{(i)} Due to the optimality of $\bbb{x}^{t+1}$, we have: $F(\bbb{x}^{t+1}) + \tfrac{\theta}{2} \|\bbb{x}^{t+1}-\bbb{x}^t\|_2^2 \leq F(\bbb{u})+ \frac{\theta}{2}\|\bbb{u}-\bbb{x}^t\|_2^2$ for all $\bbb{u}_N=\bbb{x}_N^t$. Letting $\bbb{u}=\bbb{x}^t$, we obtain the sufficient decrease condition in (\ref{eq:suff:dec}).

%

\noi Taking the expectation of $B$ for the sufficient descent inequality, we have $\textstyle\E[F(\bbb{x}^{t+1}) ] \leq F(\bbb{x}^t) - \E[\tfrac{\theta}{2}\|\bbb{x}^{t+1}-\bbb{x}^t\| ]$. Summing this inequality over $i=0,1,2, ...,t-1$, we have: $$\textstyle \tfrac{\theta}{2} \sum_{i=0}^t \E[\| \bbb{x}^{i+1}-\bbb{x}^i \|_2^2 ] \leq F(\bbb{x}^0)- F(\bbb{x}^t).$$

\noi Using the fact that $ F(\bar{\bbb{x}}) \leq F(\bbb{x}^t)$, we obtain:
\beq \label{eq:conv:bound:xx}
&&\textstyle \min_{i=1,...,t} \E[\tfrac{\theta}{2}\| \bbb{x}^{i+1}-\bbb{x}^i \|_2^2] \nn\\
&\leq& \textstyle \tfrac{\theta}{2t} \sum_{i=0}^t \E[\| \bbb{x}^{i+1}-\bbb{x}^i \|_2^2] \leq    \tfrac{F(\bbb{x}^0) - F( \bar{\bbb{x}})}{t}
\eeq
\noi Therefore, we have $\lim_{t\rightarrow \infty}\E[\|\bbb{x}^{t+1} - \bbb{x}^t\|] = 0$.

\bbb{(ii)} We assume that the stationary point is not a block-$k$ stationary point. In expectation there exists a block of coordinates $B$ such that $\bbb{x}^t \notin\arg\min_{\bbb{z}}~\P(\bbb{z};\bbb{x}^t,B)$ for some $B$, where $\P(\cdot)$ is defined in Definition \ref{def:block:k}. However, according to the fact that $\bbb{x}^t=\bbb{x}^{t+1}$ and step S2 in Algorithm \ref{algo:main}, we have $\bbb{x}^{t+1} \in\arg\min_{\bbb{z}}~\P(\bbb{z};\bbb{x}^t,B)$. Hence, we have $\bbb{x}^t_{B}\neq \bbb{x}^{t+1}_{B}$. This contradicts with the fact that $\bbb{x}^t=\bbb{x}^{t+1}$ as $t\rightarrow \infty$. We conclude that $\bbb{x}^{t}$ converges to the block-$k$ stationary point.

\bbb{(iii)} We observe that when the current solution $\bbb{x}^t$ changes, we have $\|\bbb{x}^{t+1}-\bbb{x}^t\| \geq \sqrt{2}$ if $\bbb{x}^t\neq\bbb{x}^{t+1}$. Noticing that there are $C_n^k$ possible combinations of choice for the working set of size $k$, we obtain: $\E[\|\bbb{x}_{\bbb{B}}\|_2^2] = \frac{k}{n}\|\bbb{x}\|_2^2,~\forall~\bbb{x}$. Thus, we have
\beq
\E[\|(\bbb{x}^{t+1}-\bbb{x}^t)_B\|_2^2] = \frac{k}{n}\| \bbb{x}^{t+1}-\bbb{x}^t\|_2^2 \geq \frac{2k}{n}. \nn
\eeq
\noi From (\ref{eq:conv:bound:xx}), we obtain: $\frac{F(\bbb{x}^0) - F(\bar{\bbb{x}})}{t\theta}\geq \frac{{2}k}{ n}$. Therefore, the number of iterations is upper bounded by $t \leq \frac{n [F(\bbb{x}^0) - F(\bar{\bbb{x}})]}{2k \theta}$ times in expectation.

\bbb{(iv)} Note that Algorithm \ref{algo:main} solves problem (\ref{eq:subprob}) in every iteration. Using Proposition \ref{proposition:hierarchy}, we have that the solution $\bbb{x}^{t+1}_B$ is also a $L$-stationary point. Therefore, we have $|\bbb{x}^{t+1}|_i \geq \min(\rho,\sqrt{{2\lambda}/{(\theta + L)}})$ for all $\bbb{x}^{t+1}_i\neq 0$. Taking the initial point of $\bbb{x}$ for consideration, we have that:
$$|\bbb{x}_i^{t+1}|\geq \min(\rho,|\bbb{x}_i^0|,\sqrt{2\lambda/(\theta + L)}),~\forall~ i=1,~2,...,n.$$
\noi Therefore, we have the following results: $\|\bbb{x}^{t+1}-\bbb{x}^{t}\|_2 \geq \delta$. Taking the expectation of $B$, we have the following results: $\E[\|(\bbb{x}^{t+1}-\bbb{x}^{t})_B\|_2^2] = \frac{k}{n}\|\bbb{x}^{t+1}-\bbb{x}^{t}\|_2^2 \geq \frac{k}{n} \delta^2$. Every time the index set of $\bbb{x}$ is changed, the objective value is decreased at least by $\E[\tfrac{\theta}{2}\|\bbb{x}^{t+1}-\bbb{x}^{t}\|^2]\geq \frac{k \theta\delta^2}{2n} \triangleq D$. Moreover, we obtain: $\frac{[2F(\bbb{x}^0) - 2F(\bar{\bbb{x}})]}{t\theta}  \geq \frac{\delta^2 k}{n}$ from (\ref{eq:conv:bound:xx}). Therefore, the number of iterations is upper bounded by ${\bar{J}}$.

\end{proof}
\end{theorem}

\noi \textbf{Remarks:} Coordinate descent may cycle indefinitely if each minimization step contains multiple solutions \cite{Powell1973} . The introduction of the strongly convex parameter $\theta>0$ is necessary for our nonconvex optimization problem since it guarantees sufficient decrease condition which is important for global convergence. Our algorithm is guaranteed to find the block-$k$ stationary point but it is in expectation.

\section{Additional Discussions} \label{sect:dis}

This section provides additional discussions for the proposed method.

\subsection{When the objective function is complicated}

In step $(S2)$ of the proposed algorithm, a global solution is to be found for the subproblem. When $f(\cdot)$ is simple (e.g. a quadratic function), we can find efficient and exact solutions to the subproblems. We now consider the situation when $f$ f is complicated (e.g., logistic regression, maximum entropy models). One can still find a quadratic majorizer $Q(\bbb{x},\bbb{z})$ for the convex function $f(\bbb{x})$ with
$$f(\bbb{x}) \leq Q(\bbb{x},\bbb{z}) \triangleq f(\bbb{z}) + (\bbb{x}-\bbb{z})^T\nabla f(\bbb{z})  + \tfrac{1}{2}(\bbb{x}-\bbb{z})^T\bbb{M}(\bbb{z})(\bbb{x}-\bbb{z}),~\forall~\bbb{z},~\bbb{x},~\bbb{M}(\bbb{z})\succ \nabla f^2(\bbb{z}).$$
\noi By minimizing the upper bound of $f(\bbb{x})$ (i.e., the quadratic surrogate function) at the current estimate $\bbb{x}^t$, i.e.,
$$\bbb{x}^{t+1}\Leftarrow \arg\min_{\bbb{x}}~Q(\bbb{x},\bbb{x}^{t})+h(\bbb{x}),$$
 \noi we can drive the objective downward until a stationary point is reached. We will obtain a stationary point $\ddot{\bbb{x}}$ satisfying:
 $$\ddot{\bbb{x}} = \arg \min_{\bbb{z}}~h(\bbb{z}) + f(\ddot{\bbb{x}}) + (\bbb{z}-\ddot{\bbb{x}})^T\nabla f(\ddot{\bbb{x}})  + \frac{1}{2} (\bbb{z}-\ddot{\bbb{x}})^T  \bbb{M}(\ddot{\bbb{x}}) (\bbb{z}-\ddot{\bbb{x}}),~s.t.~\ddot{\bbb{x}}_N=\bbb{z}_N$$
  \noi for all $N$. However, this is weaker than the block-$k$ stationary point $\bar{\bbb{x}}$: $\bar{\bbb{x}} = \arg \min_{\bbb{z}}~h(\bbb{z}) + f({\bbb{z}}),~s.t.~\bar{\bbb{x}}_N=\bbb{z}_N$ for all $N$.

\subsection{Practical computational efficiency}

Block coordinate descent is shown to be very efficient for solving convex problems (e.g., support vector machines \cite{chang2008coordinate,hsieh2008dual}, LASSO problems \cite{tseng2009coordinate}, nonnegative matrix factorization \cite{hsieh2011fast}). The main difference of our block coordinate descent from existing ones is that our method needs to solve a small-sized NP-hard subproblem globally which takes subexponential
time $\mathcal{O}(2^k)$. As a result, our algorithm finds a block-$k$ approximation solution
for the original NP-hard problem within $\mathcal{O}(2^k)$ time. This is expected since otherwise we have proven NP=P. When $k$ is large, it is hard to enumerate the full binary tree since the subproblem is equally NP-hard. However, $k$ is relatively small in practice (e.g., 2 to 20). In addition, real-world applications often have some special (e.g. unbalanced, sparse, local) structure and block-$k$ stationary point could also be the global stationary point (refer to Table \ref{tab:optimality} of this paper).


\begin{figure*} [!t]
\captionsetup{singlelinecheck = on, format= hang, justification=justified, font=footnotesize, labelsep=space}
\centering

      \begin{subfigure}{\fourfigwid}\includegraphics[width=\objimgwid,height=\objimghei]{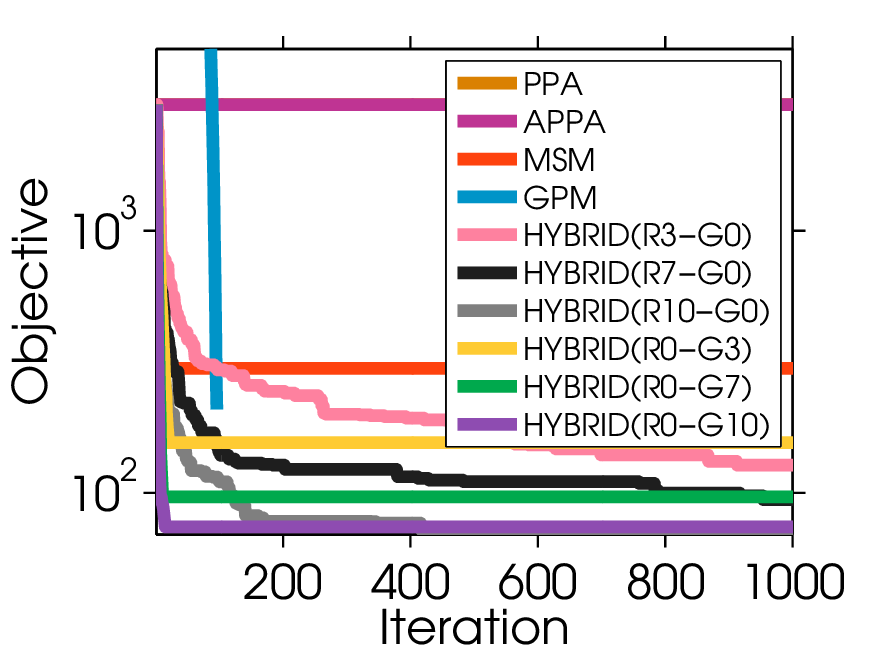}\vspace{-6pt}\caption{ $\sigma=0.001,~h_{\text{binary}}$}\end{subfigure}\ghs
      \begin{subfigure}{\fourfigwid}\includegraphics[width=\objimgwid,height=\objimghei]{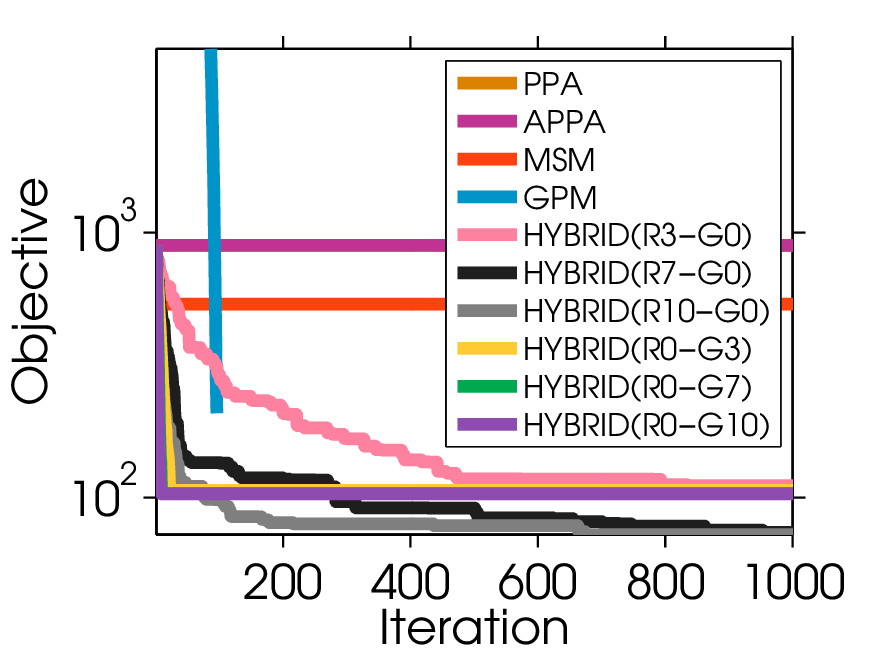}\vspace{-6pt} \caption{ $\sigma=0.01,~h_{\text{binary}}$}\end{subfigure}\ghs
      \begin{subfigure}{\fourfigwid}\includegraphics[width=\objimgwid,height=\objimghei]{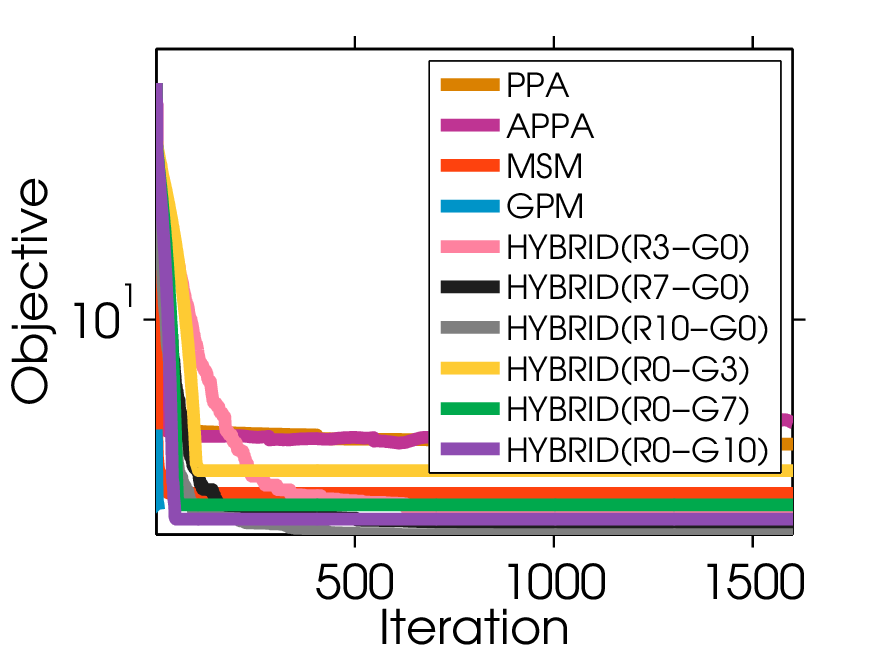}\vspace{-6pt}\caption{ $\sigma=0.001,~h_{\text{sparse}}$}\end{subfigure}\ghs
      \begin{subfigure}{\fourfigwid}\includegraphics[width=\objimgwid,height=\objimghei]{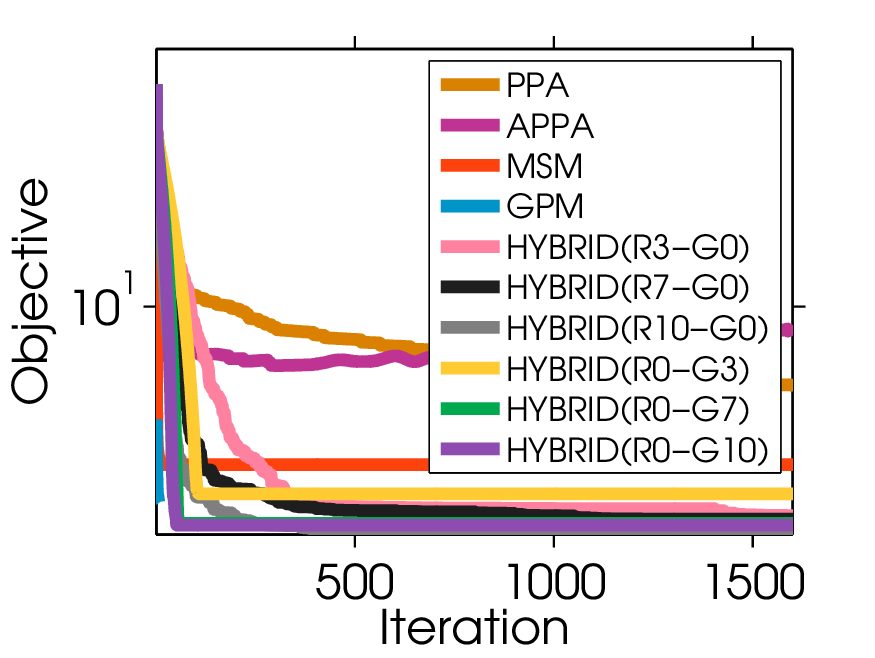}\vspace{-6pt} \caption{ $\sigma=0.01,~h_{\text{sparse}}$}\end{subfigure}

      \centering
\caption{Convergence behavior for solving the sparse regularized least squares / binary constrained least squares optimization problem with different initializations. Denoting $\tilde{\bbb{o}}$ as an arbitrary standard Gaussian random matrix of suitable size, we consider the following strategy for different initiations $\bbb{x}$. Figure a and b: $\bbb{x}= sign(\sigma \times\tilde{\bbb{o}})$. Figure c and d: $\bbb{x}=\sigma \times \tilde{\bbb{o}}$.}
\label{fig:binary:c:sparse:r}
\end{figure*}

\section{Experimental Validation} \label{sect:exp}
In this section, we demonstrate the effectiveness of our proposed algorithm on three discrete optimization tasks, namely sparse regularized least square problem, binary constrained least square problem, and sparsity constrained least square problem. We use \textbf{HYBRID (R$i$-G$j$)} to denote our hybrid method along with selecting $i$ coordinate using \bbb{R}andom strategy and $j$ coordinates using the \bbb{G}reedy strategy. Since these two strategies may select the same coordinates, the working set at most contains $i+j$ coordinates. We set $\theta=10^{-3}$ for HYBRID. We keep a record of the relative changes of the objective function values by $r_t = (f(\bbb{x}^t)-f(\bbb{x}^{t+1}))/f(\bbb{x}^t)$. We let our algorithms run up to $T$ iterations and stop them at iteration $t<T$ if $\text{mean}([{r}_{t-\text{min}(t,M)+1},~{r}_{t-min(t,M)+2},...,{r}_t]) \leq \epsilon$. The defaults values of $\epsilon$, $M$, and $T$ are $10^{-5}$, $50$ and $1000$, respectively. All codes were implemented in Matlab on an Intel 3.20GHz CPU with 8 GB RAM. 

We also apply our method to solve dense subgraph discovery problem \cite{ravi1994heuristic,feige2001dense,yuan2013truncated,wu2016ellp,YuanG17aaai} and our experiments have shown that our method achieves state-of-the-art performance. Due to space limitation, we place them into the \bbb{supplementary material}.

\subsection{Binary Constrained / Sparse Regularized Least Squares Problem}

Given a design matrix $\bbb{A} \in \mathbb{R}^{m\times n}$ and an observation vector $\bbb{b} \in \mathbb{R}^m$, sparse regularized / binary constrained least square problem is to solve the following optimization problem:
\beq
\textstyle \min_{\bbb{x}\in\{-1,1\}^n}\frac{1}{2}\|\bbb{Ax}-\bbb{b}\|_2^2~\text{or}~\min_{\bbb{x}}\frac{1}{2}\|\bbb{Ax}-\bbb{b}\|_2^2 +  \lambda\|\bbb{x}\|_0 \nn
\eeq

We generate $\bbb{A} \in \mathbb{R}^{200\times 500}$ and $\bbb{b} \in \mathbb{R}^{200}$ from a (0-1) uniform distribution. We set $\lambda=0.1$. Note that although HYBRID uses some randomized strategies to find the working set, we can always measure the quality of the solution by computing the deterministic objective value.

\bbb{Compared Methods.} We compare the proposed method (HYBRID) with four state-of-the-art methods: \bbb{(i)} Proximal Point Algorithm (PPA) \cite{nesterov2013introductory}, \bbb{(ii)} Accelerated PPA (APPA) \cite{nesterov2013introductory,beck2009fast}, \bbb{(iii)} Matrix Splitting Method (MSM) \cite{yuan2017matrix}, and \bbb{(iv)} Greedy Pursuit Method (GPM).

\bbb{Experimental Results.} Several observations can be drawn from Figure \ref{fig:binary:c:sparse:r}. \bbb{(i)} PPA and APPA achieve nearly the same performance and they get stuck into poor local minima. \bbb{(ii)} MSM improves over PGA and APGA, and GPM consistently outperforms MSM. \bbb{(iii)} Our proposed hybrid method is more effective than MSM and GPM. In addition, we find that as the parameter $k$ becomes larger, more higher accuracy is achieved. \bbb{(iii)} HYBRID appears to be less sensitive to initialization and it converges to similar objective values when using different initializations. \bbb{(iv)} We notice that Hybrid(R0-G10) converges quickly but it generally leads to worse solution quality than HYBRID(R10-G0). Based on this observation, we consider a combined random and greedy strategies for finding the working set in our forthcoming experiments.

\begin{figure*} [!t]
\centering
      \begin{subfigure}{\fourfigwid}\includegraphics[width=\objimgwid,height=\objimghei]{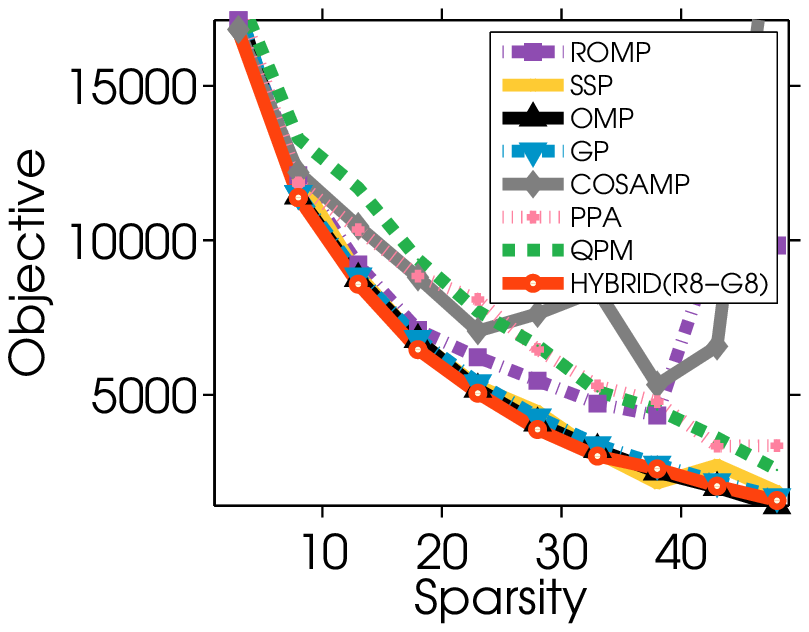}\vspace{-6pt} \caption{ $n=256$ }\end{subfigure}\ghs
      \begin{subfigure}{\fourfigwid}\includegraphics[width=\objimgwid,height=\objimghei]{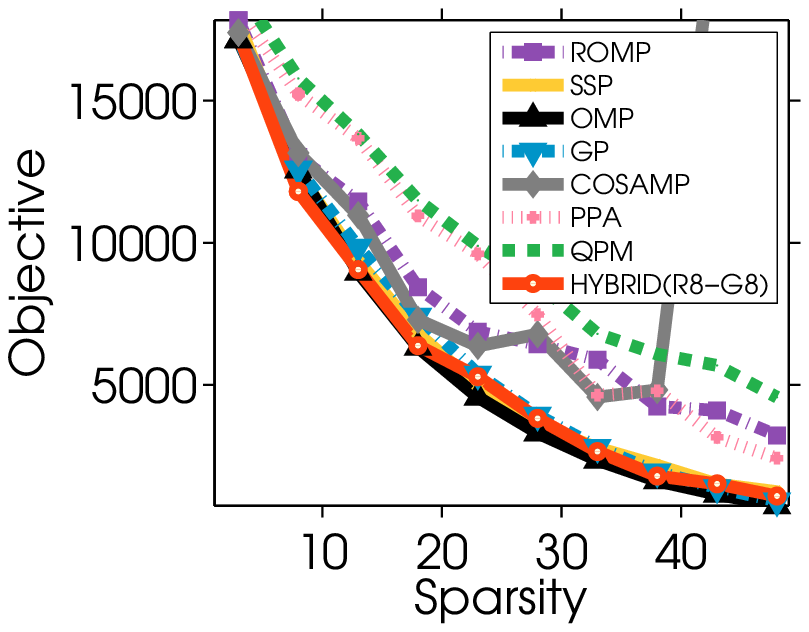}\vspace{-6pt} \caption{ $n=512$}\end{subfigure}\ghs
      \begin{subfigure}{\fourfigwid}\includegraphics[width=\objimgwid,height=\objimghei]{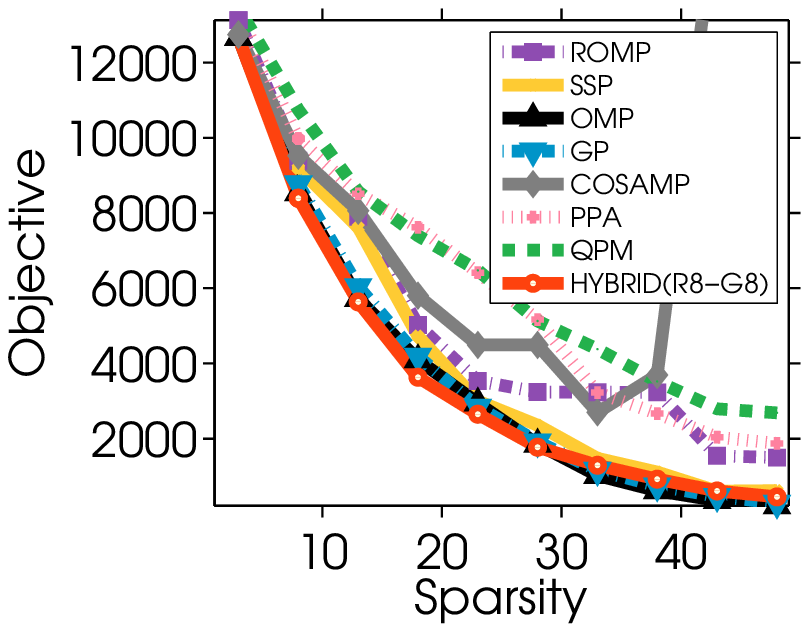}\vspace{-6pt} \caption{ $n=1024$}\end{subfigure}\ghs
      \begin{subfigure}{\fourfigwid}\includegraphics[width=\objimgwid,height=\objimghei]{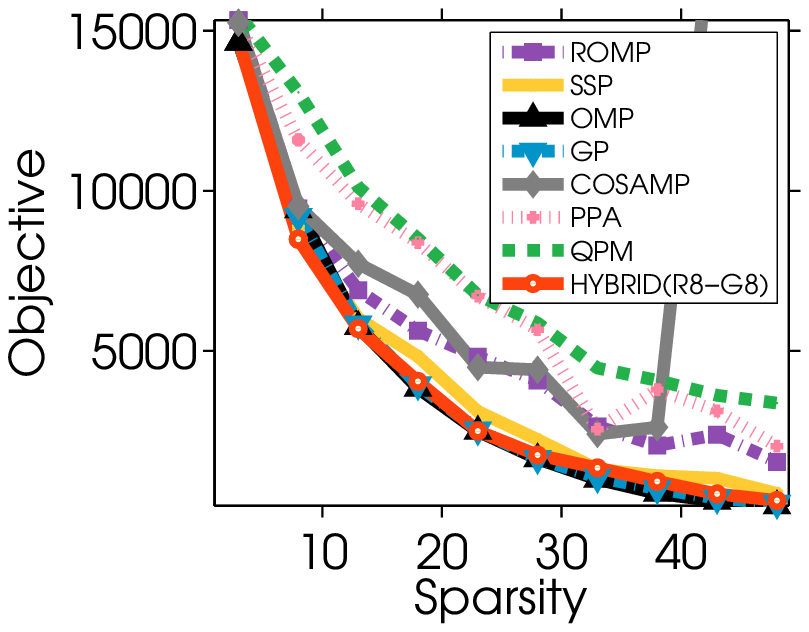}\vspace{-6pt} \caption{ $n=2048$}\end{subfigure}\\

\centering
\caption{Experimental results on sparsity constrained least squares problems on `\text{AI + bI}' with fixing $m=512$ and varying $n=\{256,~512,~1024,~2048\}$. }
\label{fig:sparse:LS:veryn:a1:b1}

\centering
      \begin{subfigure}{\fourfigwid}\includegraphics[width=\objimgwid,height=\objimghei]{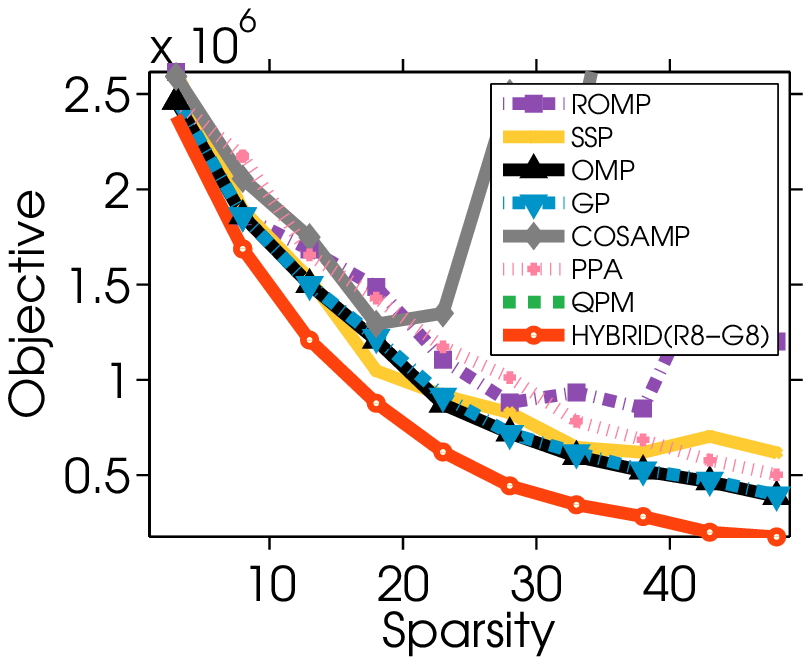}\vspace{-6pt} \caption{ $n=256$}\end{subfigure}\ghs
      \begin{subfigure}{\fourfigwid}\includegraphics[width=\objimgwid,height=\objimghei]{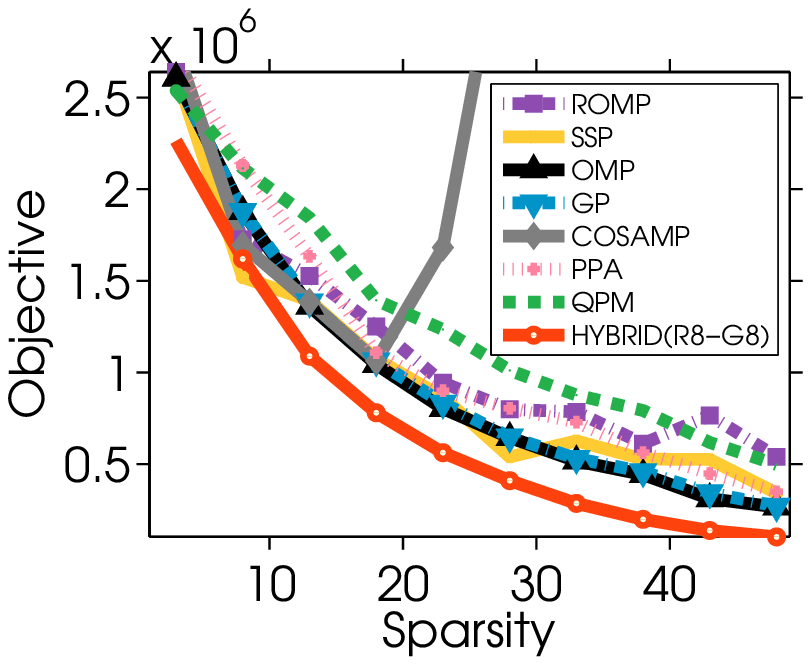}\vspace{-6pt} \caption{ $n=512$}\end{subfigure}\ghs
      \begin{subfigure}{\fourfigwid}\includegraphics[width=\objimgwid,height=\objimghei]{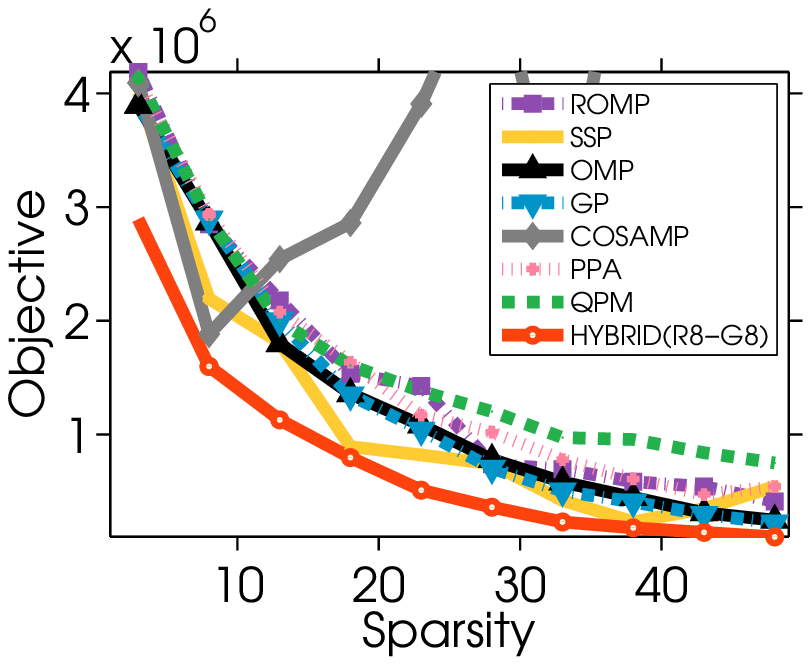}\vspace{-6pt} \caption{ $n=1024$}\end{subfigure}\ghs
      \begin{subfigure}{\fourfigwid}\includegraphics[width=\objimgwid,height=\objimghei]{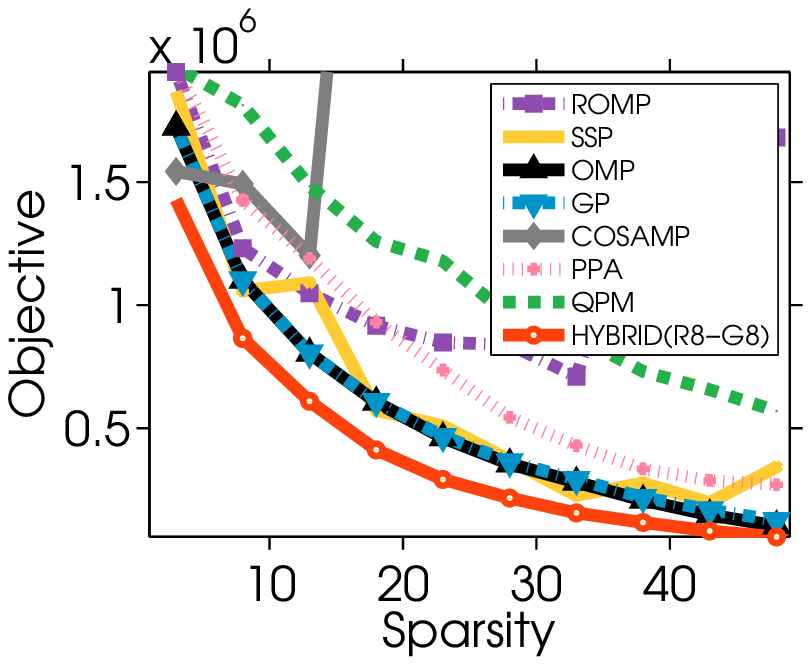}\vspace{-6pt} \caption{ $n=2048$}\end{subfigure}\\

\centering
\caption{Experimental results on sparsity constrained least squares problems on `\text{AII + bII}' with fixing $m=512$ and varying $n=\{256,~512,~1024,~2048\}$. }
\label{fig:sparse:LS:veryn:a2:b2}
\end{figure*}

\begin{figure*} [!t]
\centering

     \begin{subfigure}{\fourfigwid}\includegraphics[width=\objimgwid,height=\objimghei]{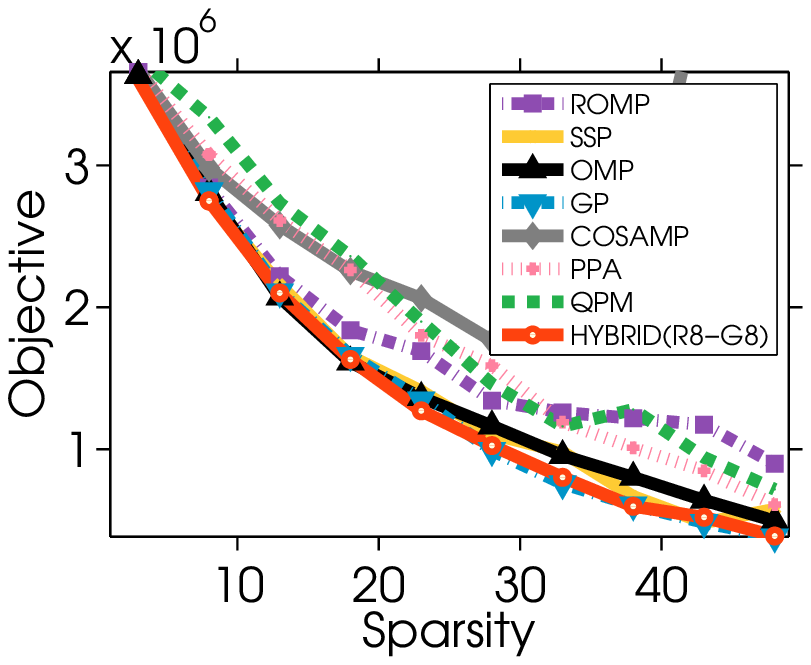}\vspace{-6pt} \caption{ $n=256$}\end{subfigure}\ghs
      \begin{subfigure}{\fourfigwid}\includegraphics[width=\objimgwid,height=\objimghei]{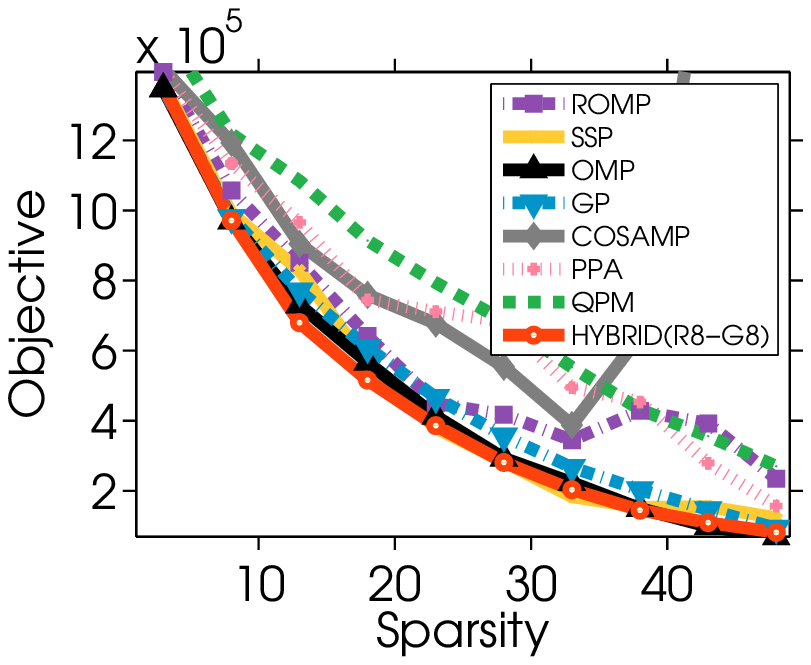}\vspace{-6pt} \caption{ $n=512$}\end{subfigure}\ghs
      \begin{subfigure}{\fourfigwid}\includegraphics[width=\objimgwid,height=\objimghei]{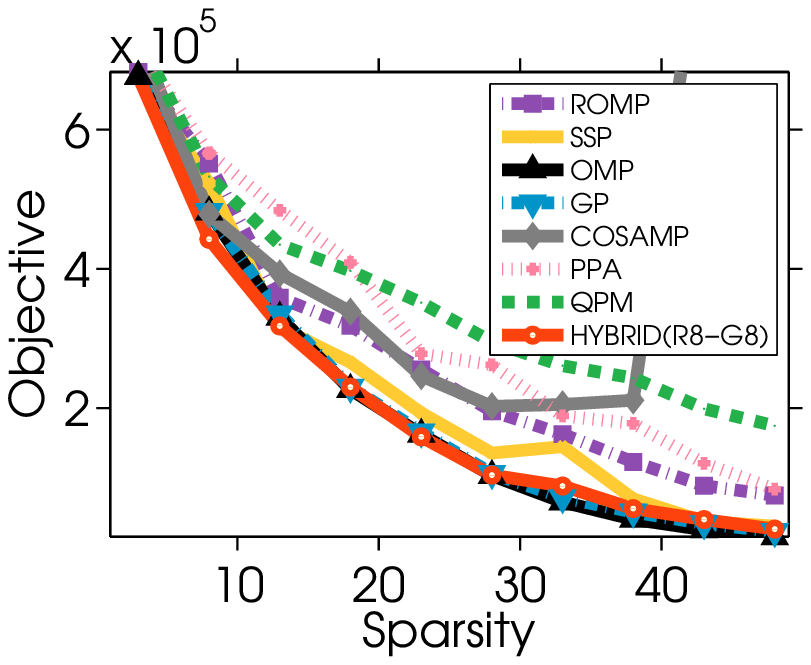}\vspace{-6pt} \caption{ $n=1024$}\end{subfigure}\ghs
      \begin{subfigure}{\fourfigwid}\includegraphics[width=\objimgwid,height=\objimghei]{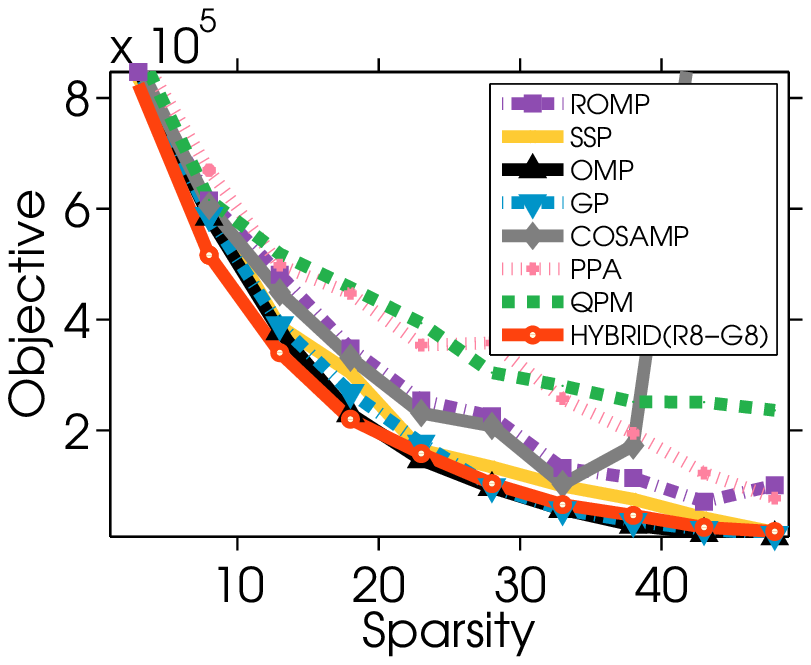}\vspace{-6pt} \caption{ $n=2048$}\end{subfigure}\\

\centering
\caption{Experimental results on sparsity constrained least squares problems on `\text{AI + bII}' with fixing $m=512$ and varying $n=\{256,~512,~1024,~2048\}$. }
\label{fig:sparse:LS:veryn:a1:b2}

\centering

      \begin{subfigure}{\fourfigwid}\includegraphics[width=\objimgwid,height=\objimghei]{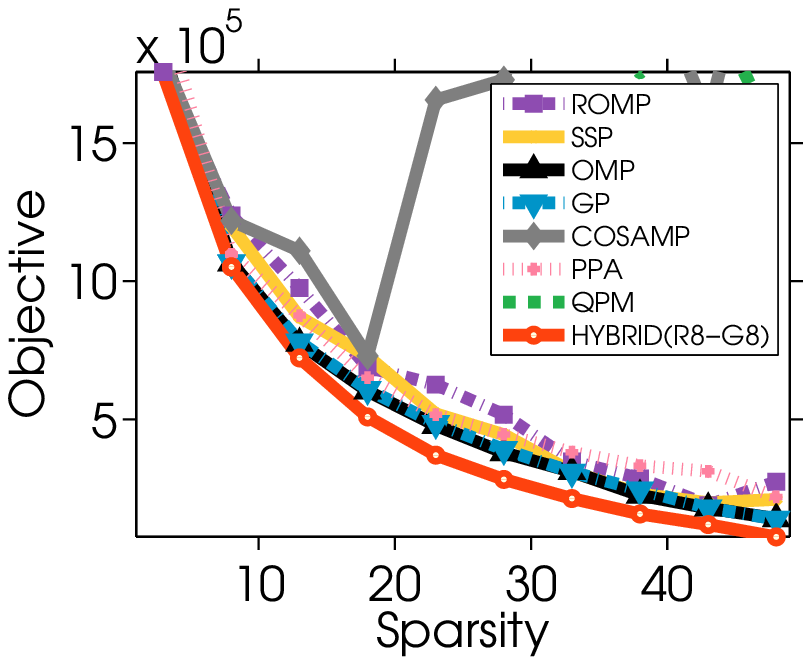}\vspace{-6pt} \caption{ $n=256$ }\end{subfigure}\ghs
      \begin{subfigure}{\fourfigwid}\includegraphics[width=\objimgwid,height=\objimghei]{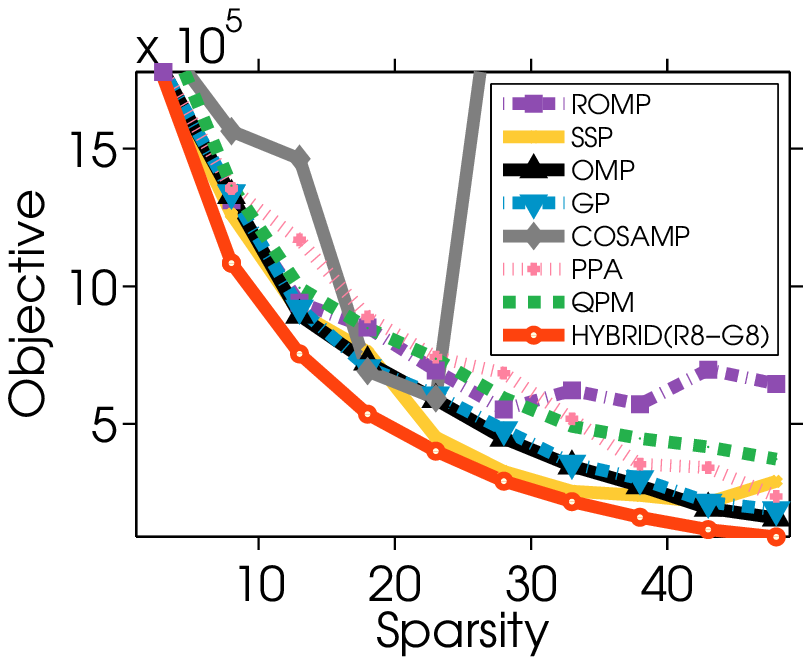}\vspace{-6pt} \caption{ $n=512$}\end{subfigure}\ghs
      \begin{subfigure}{\fourfigwid}\includegraphics[width=\objimgwid,height=\objimghei]{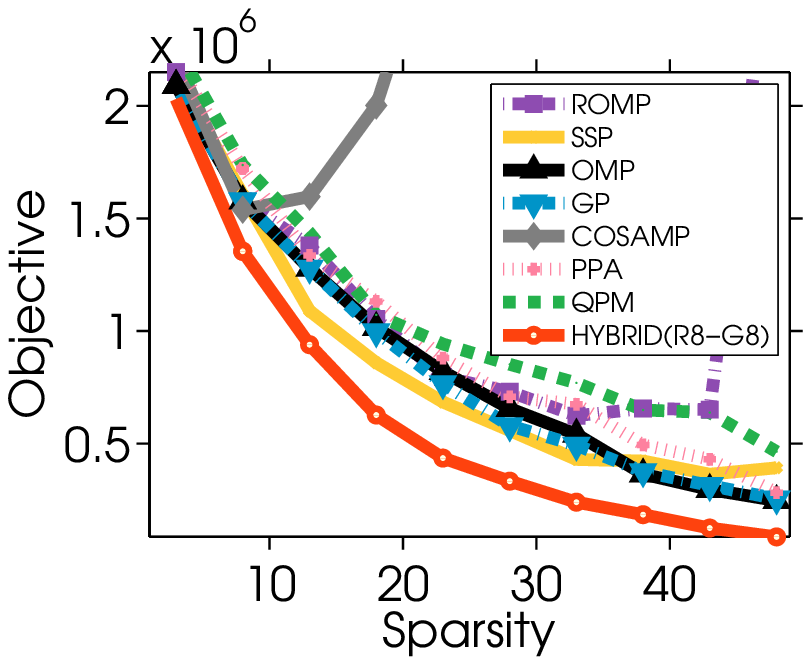}\vspace{-6pt} \caption{ $n=1024$}\end{subfigure}\ghs
      \begin{subfigure}{\fourfigwid}\includegraphics[width=\objimgwid,height=\objimghei]{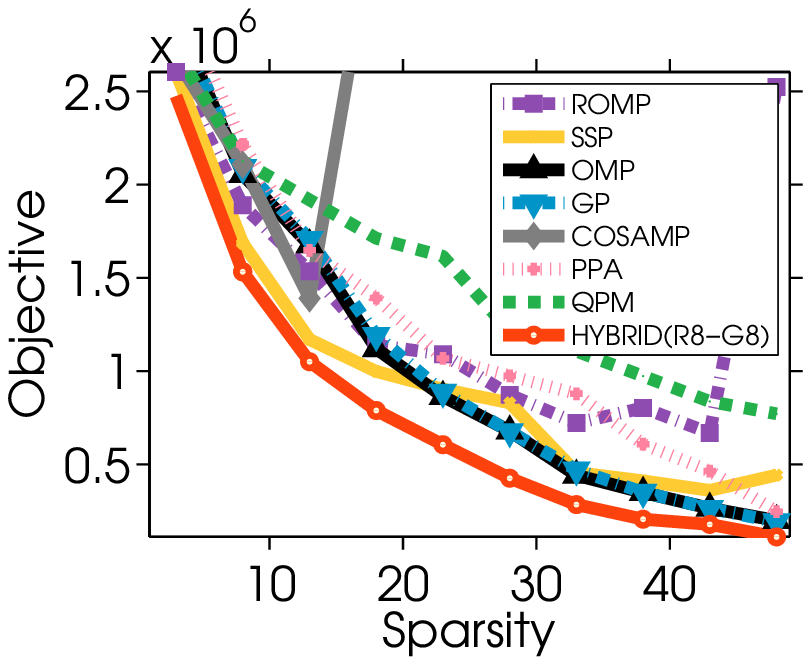}\vspace{-6pt} \caption{ $n=2048$}\end{subfigure}\\

\centering
\caption{Experimental results on sparsity constrained least squares problems on `\text{AII + bI}' with fixing $m=512$ and varying $n=\{256,~512,~1024,~2048\}$. }
\label{fig:sparse:LS:veryn:a2:b1}
\end{figure*}

\begin{figure*} [!t]
\centering
      \begin{subfigure}{\fourfigwid}\includegraphics[width=\objimgwid,height=\objimghei]{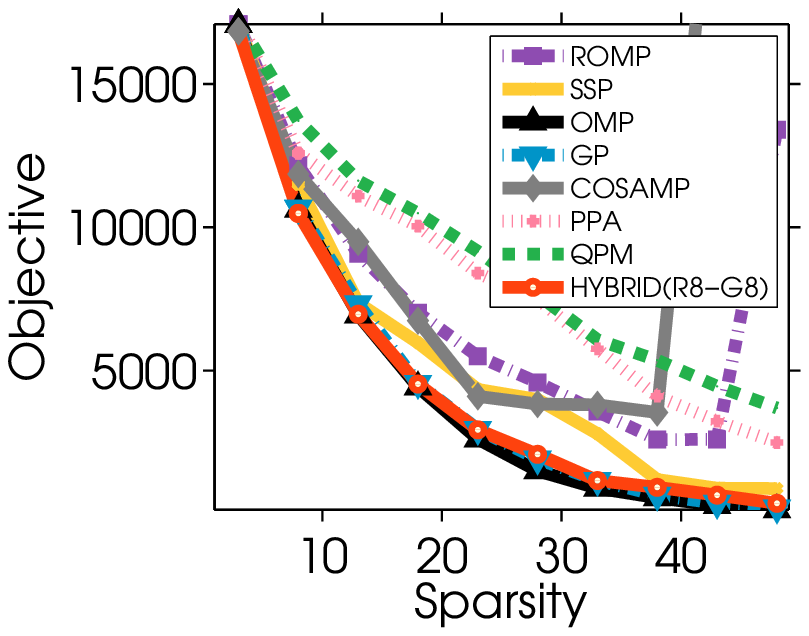}\vspace{-6pt} \caption{ $n=256$ }\end{subfigure}\ghs
      \begin{subfigure}{\fourfigwid}\includegraphics[width=\objimgwid,height=\objimghei]{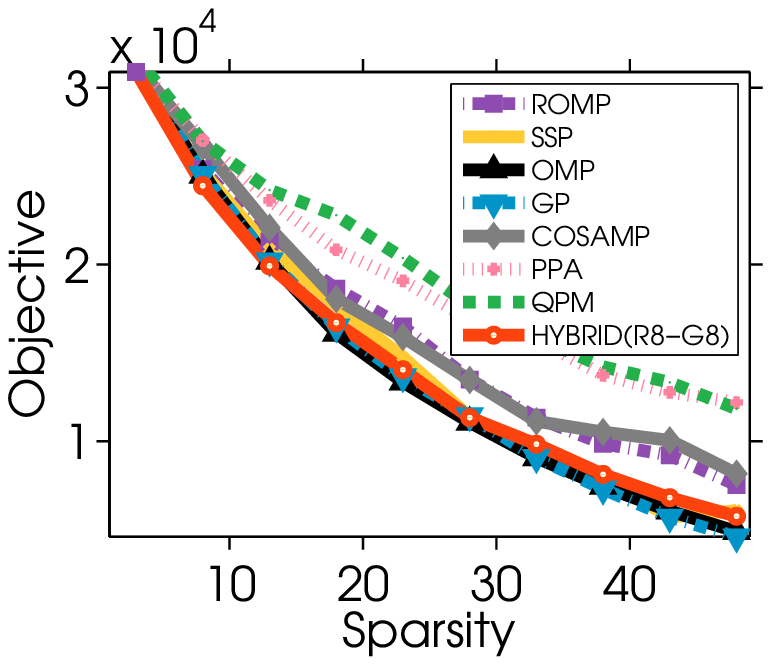}\vspace{-6pt} \caption{ $n=512$ }\end{subfigure}\ghs
      \begin{subfigure}{\fourfigwid}\includegraphics[width=\objimgwid,height=\objimghei]{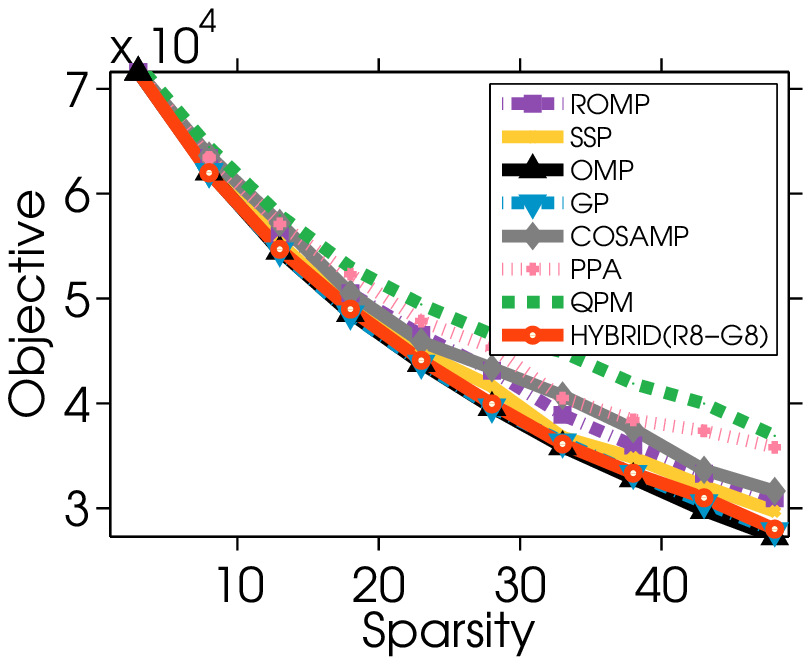}\vspace{-6pt} \caption{ $n=1024$ }\end{subfigure}\ghs
      \begin{subfigure}{\fourfigwid}\includegraphics[width=\objimgwid,height=\objimghei]{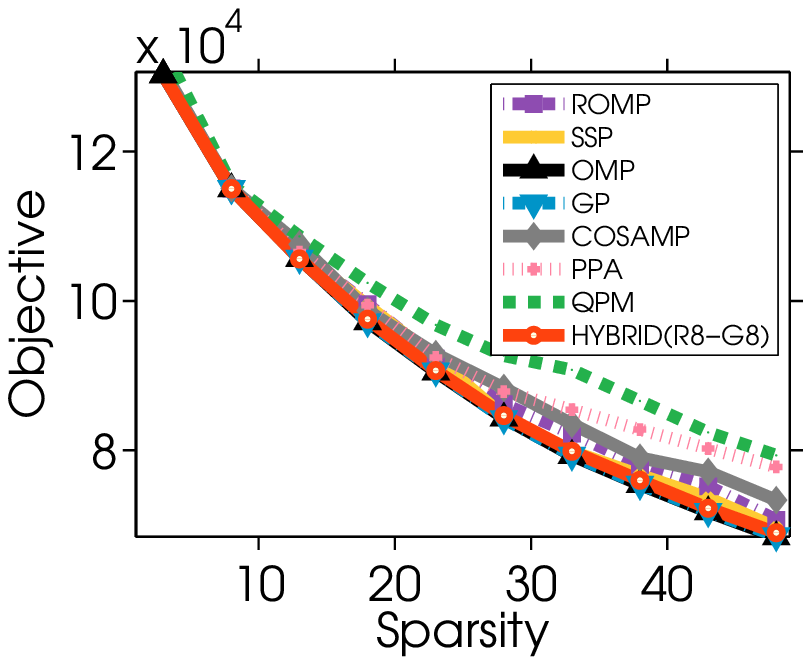}\vspace{-6pt} \caption{ $n=2048$ }\end{subfigure}\\

\centering
\caption{Experimental results on sparsity constrained least squares problems on `\text{AI + bI}' with fixing $n=2048$ and varying $m=\{128,~256,~512,~1024\}$. }
\label{fig:sparse:LS:verym:a1:b1}
%
%
\centering

      \begin{subfigure}{\fourfigwid}\includegraphics[width=\objimgwid,height=\objimghei]{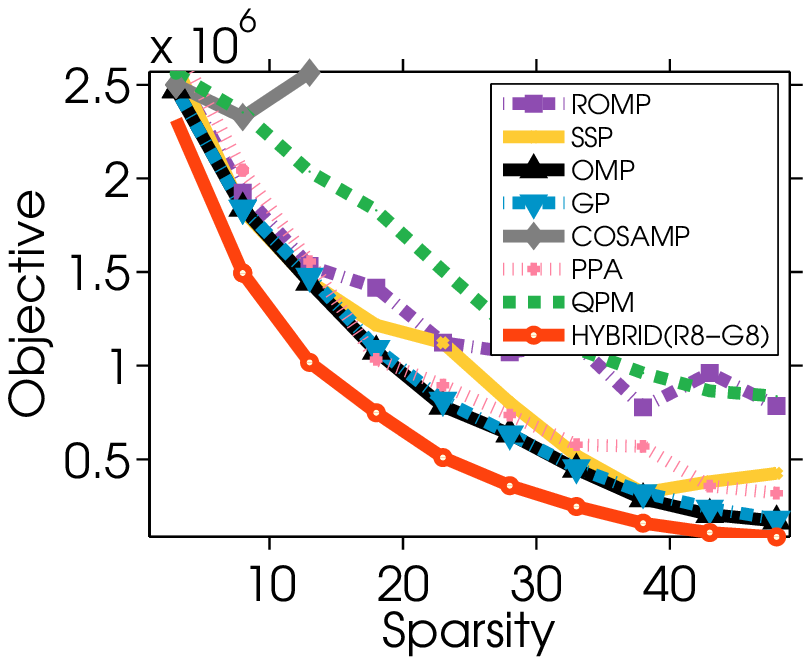}\vspace{-6pt} \caption{ $n=256$ }\end{subfigure}\ghs
      \begin{subfigure}{\fourfigwid}\includegraphics[width=\objimgwid,height=\objimghei]{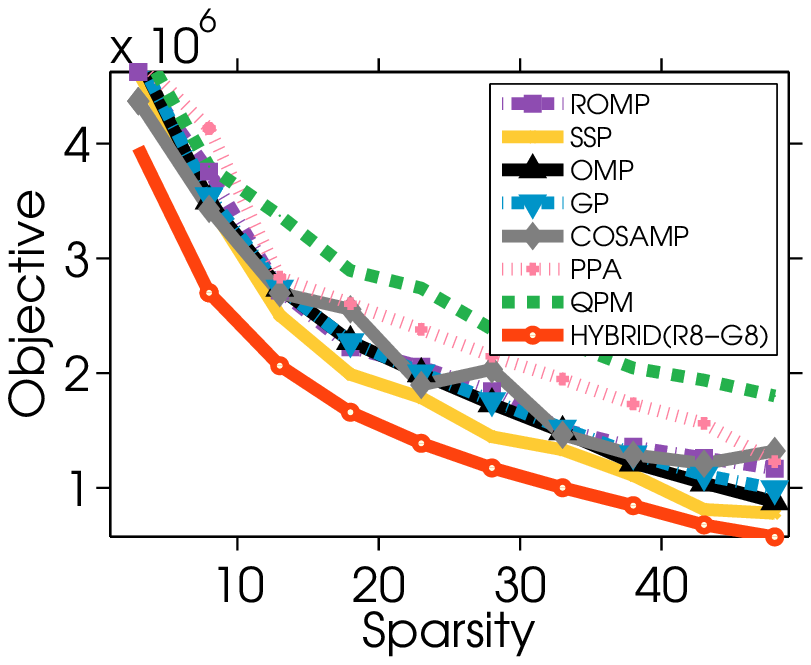}\vspace{-6pt} \caption{ $n=512$ }\end{subfigure}\ghs
      \begin{subfigure}{\fourfigwid}\includegraphics[width=\objimgwid,height=\objimghei]{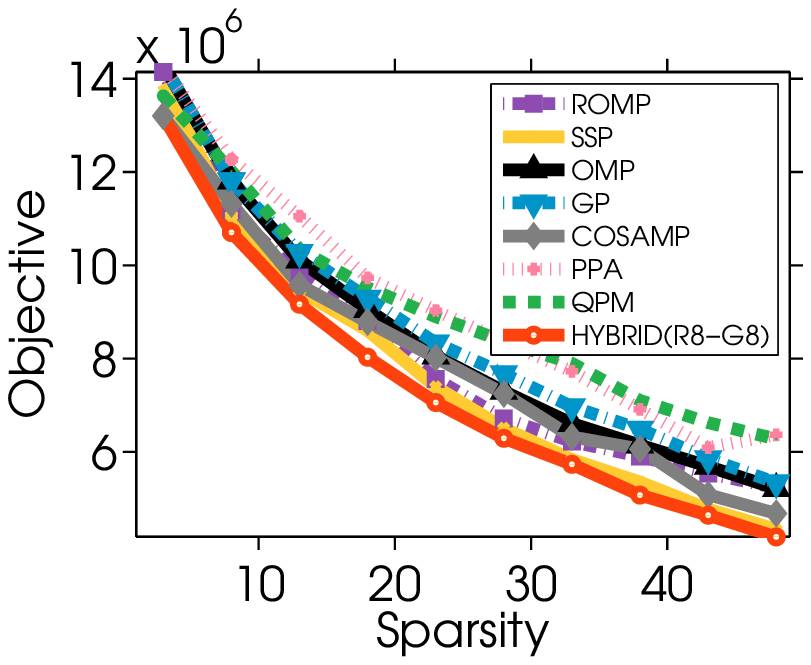}\vspace{-6pt} \caption{ $n=1024$ }\end{subfigure}\ghs
      \begin{subfigure}{\fourfigwid}\includegraphics[width=\objimgwid,height=\objimghei]{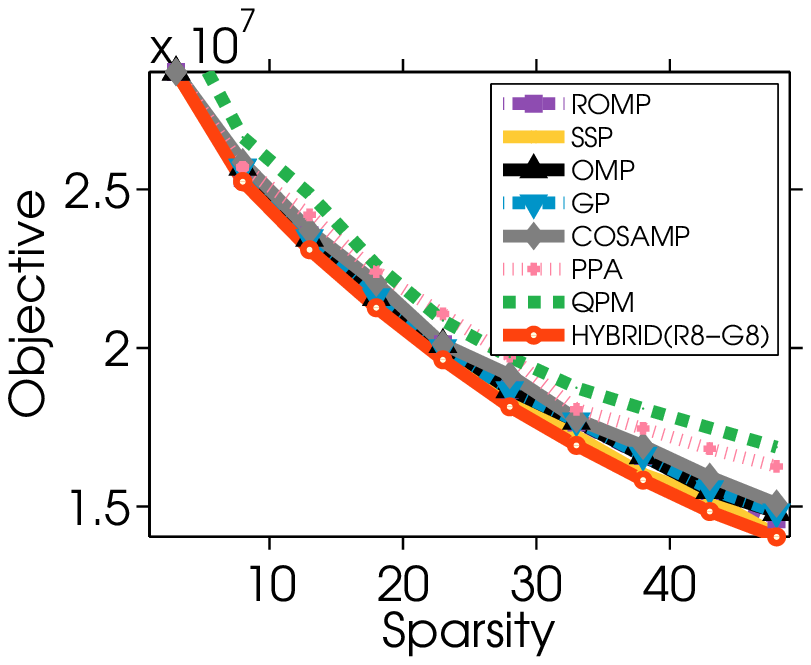}\vspace{-6pt} \caption{ $n=2048$ }\end{subfigure}\\

\centering
\caption{Experimental results on sparsity constrained least squares problems on `\text{AII + bII}' with fixing $n=2048$ and varying $m=\{128,~256,~512,~1024\}$.}
\label{fig:sparse:LS:verym:a2:b2}
\end{figure*}

\begin{figure*} [!t]
\centering

      \begin{subfigure}{\fourfigwid}\includegraphics[width=\objimgwid,height=\objimghei]{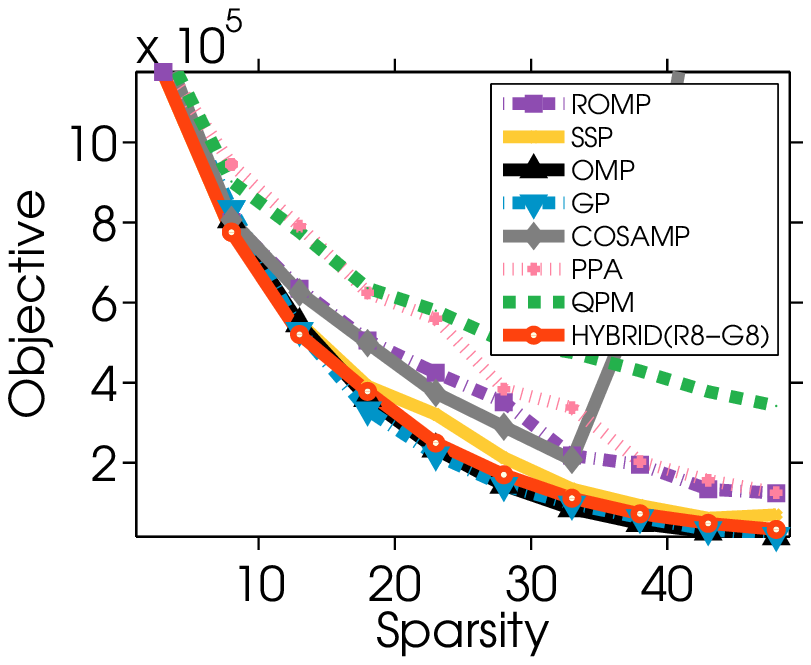}\vspace{-6pt} \caption{ $n=256$ }\end{subfigure}\ghs
      \begin{subfigure}{\fourfigwid}\includegraphics[width=\objimgwid,height=\objimghei]{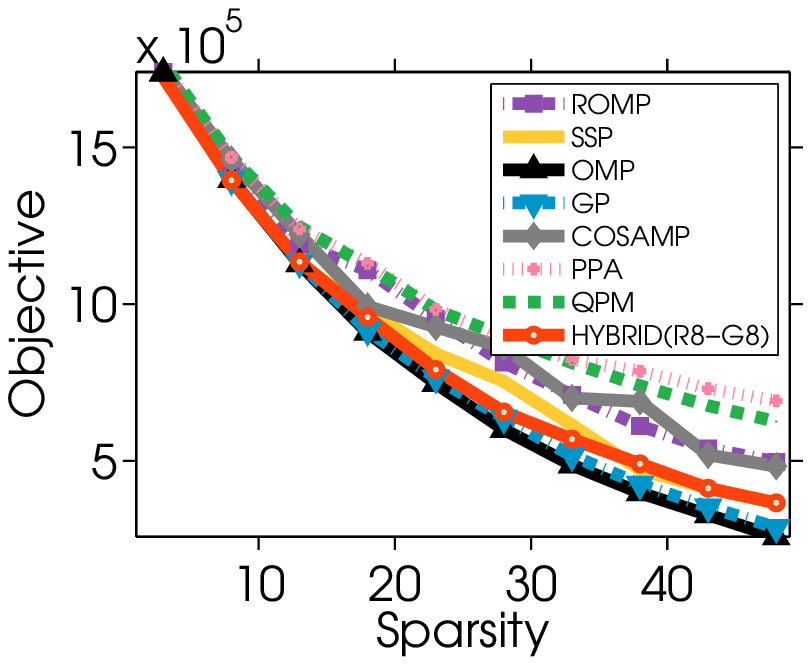}\vspace{-6pt} \caption{ $n=512$ }\end{subfigure}\ghs
      \begin{subfigure}{\fourfigwid}\includegraphics[width=\objimgwid,height=\objimghei]{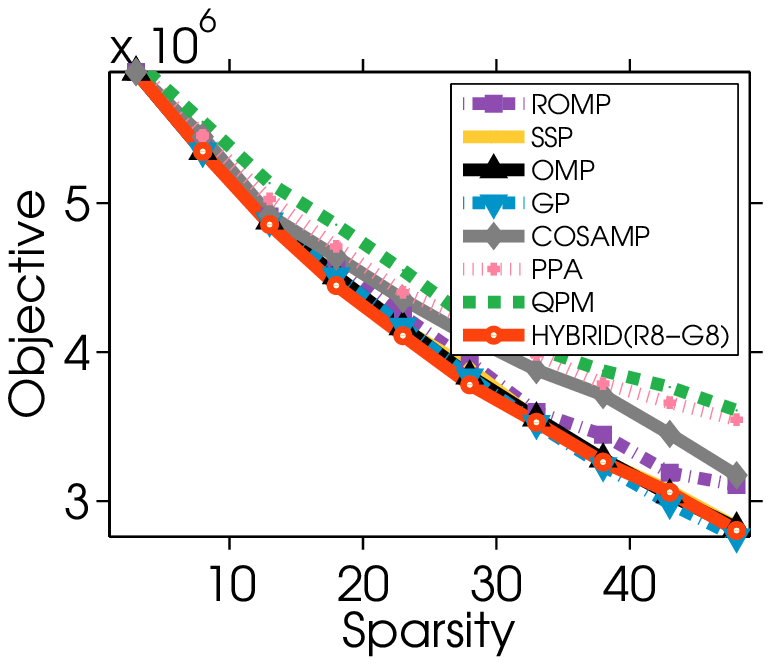}\vspace{-6pt} \caption{ $n=1024$ }\end{subfigure}\ghs
      \begin{subfigure}{\fourfigwid}\includegraphics[width=\objimgwid,height=\objimghei]{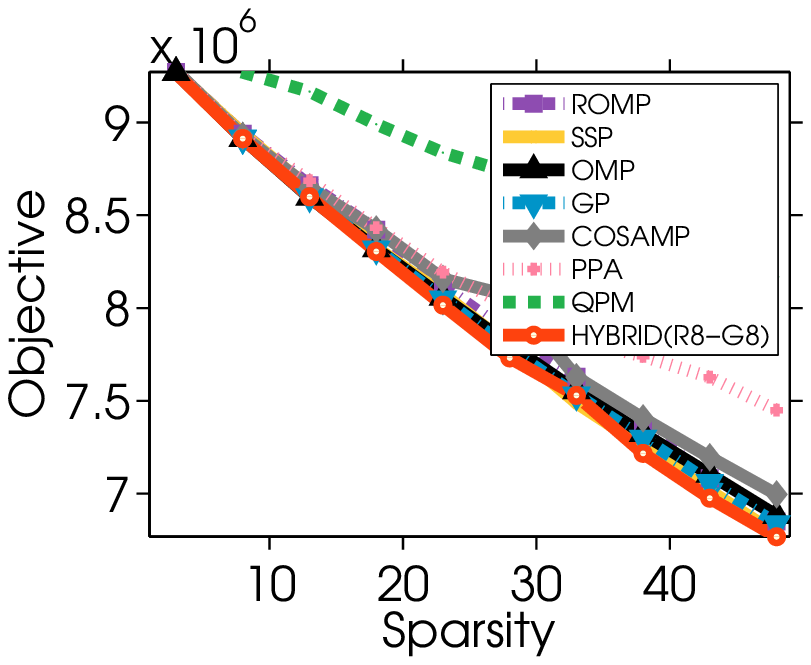}\vspace{-6pt} \caption{ $n=2048$ }\end{subfigure}\\

\centering
\caption{Experimental results on sparsity constrained least squares problems on `\text{AI + bII}' with fixing $n=2048$ and varying $m=\{128,~256,~512,~1024\}$.  }
\label{fig:sparse:LS:verym:a1:b2}
%
%
%
\centering

      \begin{subfigure}{\fourfigwid}\includegraphics[width=\objimgwid,height=\objimghei]{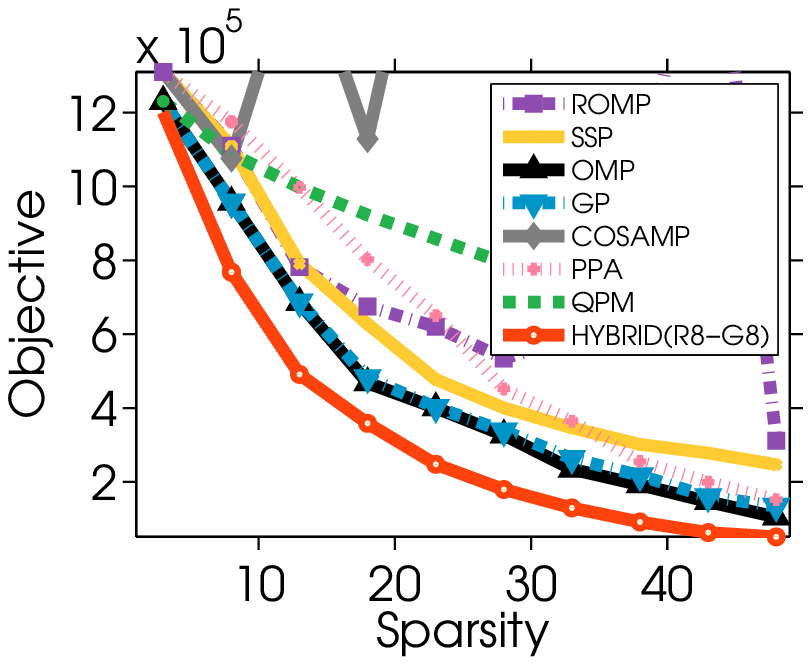}\vspace{-6pt} \caption{ $n=256$ }\end{subfigure}\ghs
      \begin{subfigure}{\fourfigwid}\includegraphics[width=\objimgwid,height=\objimghei]{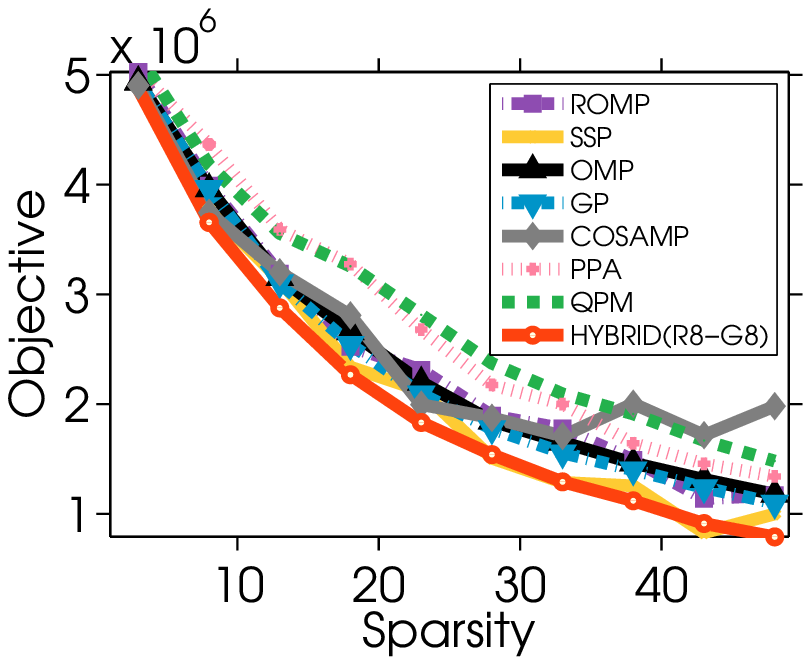}\vspace{-6pt} \caption{ $n=512$ }\end{subfigure}\ghs
      \begin{subfigure}{\fourfigwid}\includegraphics[width=\objimgwid,height=\objimghei]{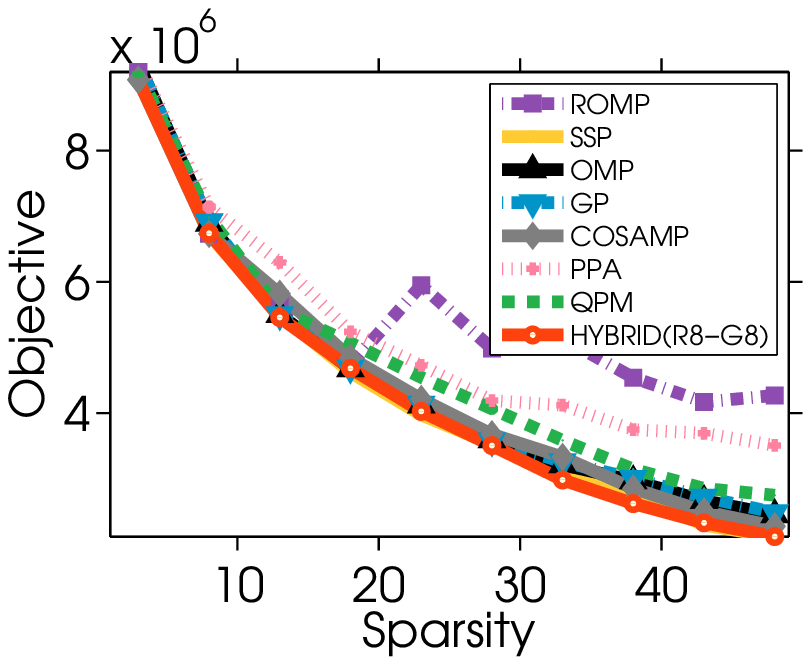}\vspace{-6pt} \caption{ $n=1024$ }\end{subfigure}\ghs
      \begin{subfigure}{\fourfigwid}\includegraphics[width=\objimgwid,height=\objimghei]{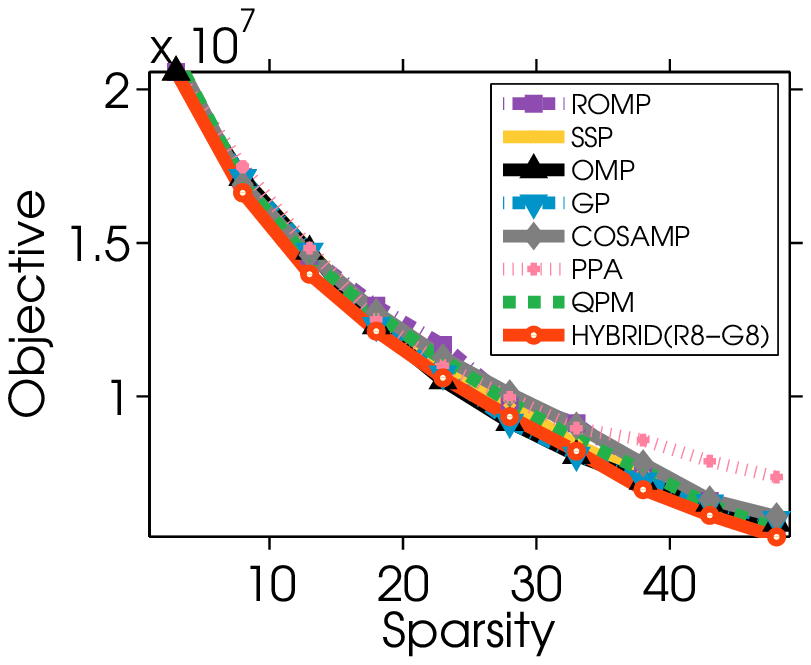}\vspace{-6pt} \caption{ $n=2048$ }\end{subfigure}\\

\centering
\caption{Experimental results on sparsity constrained least squares problems on `\text{AII + bI}' with fixing $n=2048$ and varying $m=\{128,~256,~512,~1024\}$.  }
\label{fig:sparse:LS:verym:a2:b1}
\end{figure*}

\subsection{Sparsity Constrained Least Squares Problem}

We consider the following sparsity constrained least squares problem:
\beq
\textstyle  \min_{\bbb{x}}~\frac{1}{2}\|\bbb{Ax}-\bbb{b}\|_2^2,~s.t.~\|\bbb{x}\|_0 \leq s,
\eeq
\noi In our experiments, to generate the sparse original signal $\ddot{\bbb{x}} \in \mathbb{R}^n$, we select a support set of size $100$ uniformly at random and set them to arbitrary number sampled from standard Gaussian distribution. In order to verify the robustness of the comparing methods, we generate the design matrix $\bbb{A}$ and the noise vector $\bbb{o}\in\mathbb{R}^m$ with and without outliers, as follows:
\beq
\text{AI:} ~\bbb{A} = \text{randn}(m,n),~\text{AII:} ~\bbb{A} = \P(\text{randn}(m,n));~~~~~~\nn\\
\text{bI:} ~ \bbb{o} = 10\times\text{randn}(m,1),~\text{bII:} ~\bbb{o} = \P(10\times\text{randn}(m,1)).\nn
\eeq
\noi Here, $\text{randn}(m,n)$ is a function that returns a standard Gaussian random matrix of size $m\times n$, $\P(\bbb{X})\in\mathbb{R}^{m\times p}$ denotes a noisy version of $\bbb{X}\in\mathbb{R}^{m\times p}$ where $2\%$ of the entries of $\bbb{X}$ are corrupted uniformly by scaling the original values by 100 times \footnote{Matlab script: I = randperm(m*p,round(0.02*m*p)); X(I) = X(I)*100.}. The observation vector is generated via $\bbb{b}=\bbb{A}\ddot{\bbb{x}}+\bbb{o}$. Note that the Hessian matrix can be ill-conditioned for the `AII' type design matrix. We vary $n$ from $\{256,~512,~1024,~2048\}$ and vary $m$ from $\{128,~256,~\bbb{512},~1024\}$. Unless otherwise specified, the default parameters in bold are used. We swap the parameter $s$ over $\{3,8,13,18...50\}$.

\bbb{Compared Methods.} We compare the proposed hybrid algorithm with seven state-of-the-art sparse optimization algorithms: \bbb{(i)} Regularized Orthogonal Matching Pursuit (ROMP) \cite{needell2010signal}, \bbb{(ii)} Subspace Pursuit (SSP) \cite{dai2009subspace}, \bbb{(iii)} Orthogonal Matching Pursuit (OMP) \cite{tropp2007signal}, \bbb{(iv)} Gradient Pursuit (GP) \cite{blumensath2008gradient}, \bbb{(v)} Compressive Sampling Matched Pursuit (CoSaMP)\cite{needell2009cosamp}, \bbb{(vi)} Proximal Point Algorithm (PPA) \cite{BaoJQS16}, and \bbb{(vii)} Quadratic Penalty Method (QPM) \cite{LuZ13}. We remark that ROMP, SSP, OMP, GP and CoSaMP are greedy algorithms and their support sets need to be selected iteratively. They are non-gradient type algorithms, it is hard to incorporate these methods into other gradient-type based optimization algorithms \cite{BaoJQS16}. We use the Matlab implementation in the `sparsify' toolbox\footnote{\url{http://www.personal.soton.ac.uk/tb1m08/sparsify/sparsify.html}}. Both PPA and QPM are based on iterative hard thresholding. Since the optimal solution is expected to be sparse, we initialize the solutions of PPA, QPM and HYBRID to $10^{-7} \times \text{randn}(n,1)$ and project them to feasible solutions. The initial solution of greedy pursuit methods are initialized to zero points implicitly. We show the average results of using 3 random initial points.



\bbb{Experimental Results.} We show our experimental results on sparsity constrained least squares problem with fixing $m=512$ and varying $n=\{256,~512,~1024,~2048\}$ on different type of $\bbb{A}$ and $\bbb{b}$ (see Figure \ref{fig:sparse:LS:veryn:a1:b1}, \ref{fig:sparse:LS:veryn:a2:b2}, \ref{fig:sparse:LS:veryn:a1:b2}, \ref{fig:sparse:LS:veryn:a2:b1}). We also show our the experimental results on sparsity constrained least squares problems with fixing $n=2048$ and varying $m=\{128,~256,~512,~1024\}$ on different type of $\bbb{A}$ and $\bbb{b}$ (see Figure \ref{fig:sparse:LS:verym:a1:b1},~\ref{fig:sparse:LS:verym:a2:b2},~\ref{fig:sparse:LS:verym:a1:b2},~\ref{fig:sparse:LS:verym:a2:b1}). Several conclusions can be drawn. \bbb{(i)} PPA and QPM generally lead to the worst performance. \bbb{(ii)} ROMP and COSAMP are not stable and sometimes they present bad performance. \bbb{(iii)} SSP, OMP and GP generally present comparable performance to HYBRID when the Hessian matrix is well-conditioned (for `AI' type design matrix) but present much worse performance than HYBRID when the Hessian matrix is ill-conditioned (for `AII' type design matrix). 

%

\begin{figure*} [!t]

\label{fig:exp:ls:sparse:binary}

\centering
      \begin{subfigure}{\fourfigwid}\includegraphics[width=\objimgwid,height=\objimghei]{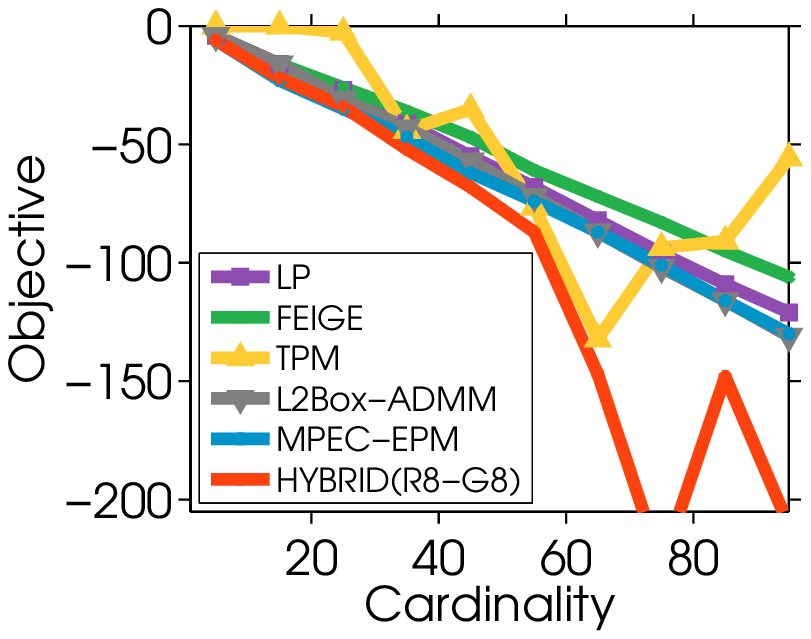}\vspace{-6pt} \caption{\footnotesize CollegeMsg}\end{subfigure}\ghs
      \begin{subfigure}{\fourfigwid}\includegraphics[width=\objimgwid,height=\objimghei]{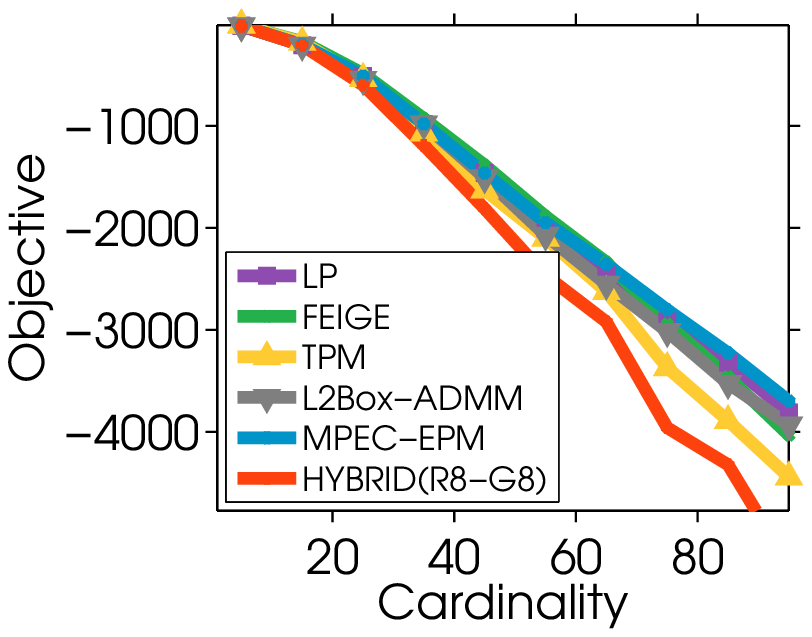}\vspace{-6pt} \caption{\footnotesize Email-Eu-Core}\end{subfigure}\ghs
      \begin{subfigure}{\fourfigwid}\includegraphics[width=\objimgwid,height=\objimghei]{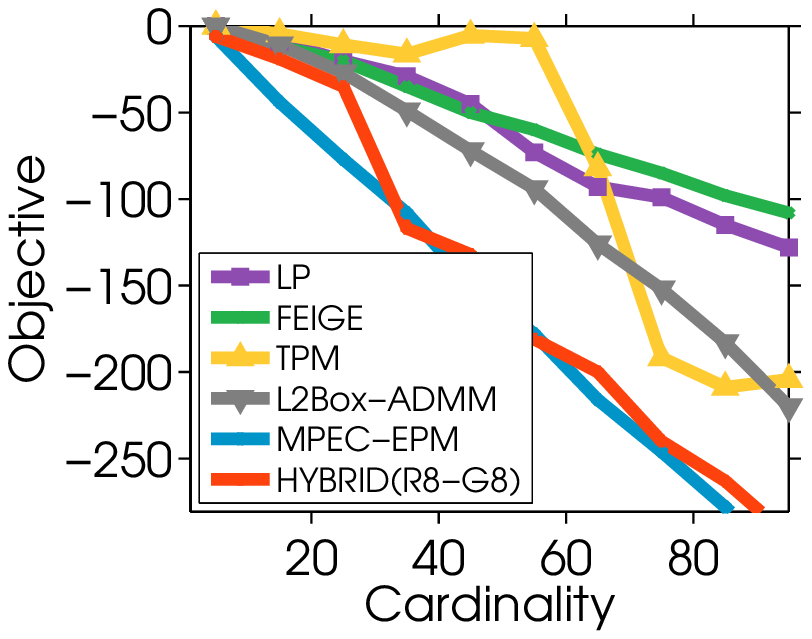}\vspace{-6pt} \caption{\footnotesize Enron}\end{subfigure}\ghs
      \begin{subfigure}{\fourfigwid}\includegraphics[width=\objimgwid,height=\objimghei]{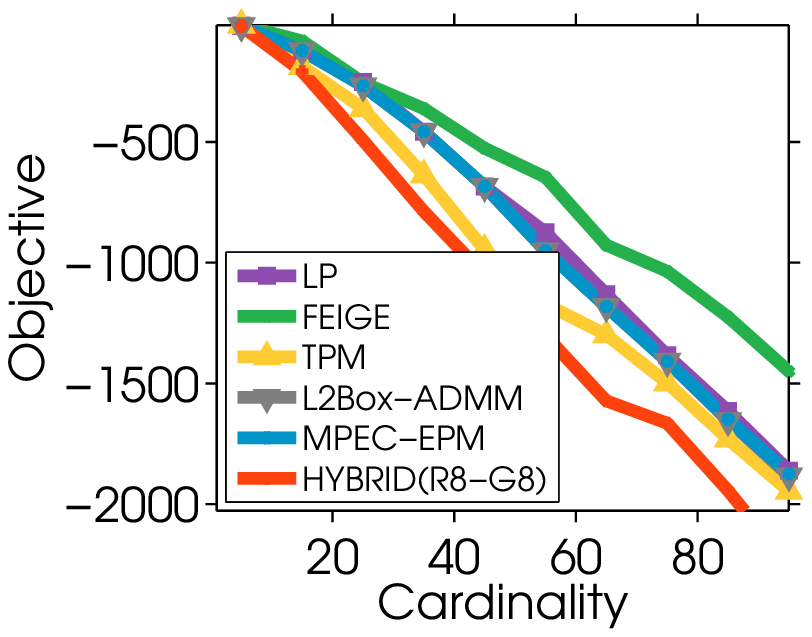}\vspace{-6pt} \caption{\footnotesize CA-CondMat}\end{subfigure}\\
      \begin{subfigure}{\fourfigwid}\includegraphics[width=\objimgwid,height=\objimghei]{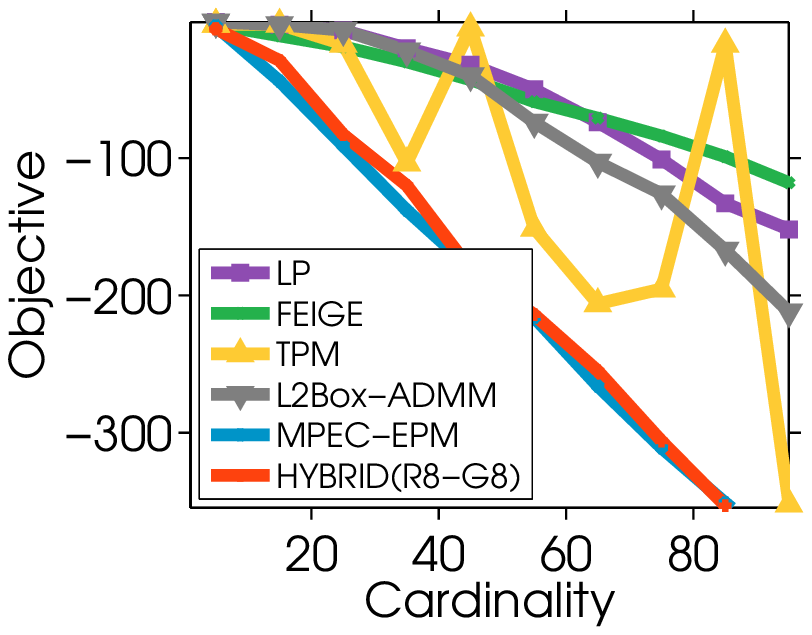}\vspace{-6pt} \caption{\footnotesize p2p-Gnutella04}\end{subfigure}\ghs
      \begin{subfigure}{\fourfigwid}\includegraphics[width=\objimgwid,height=\objimghei]{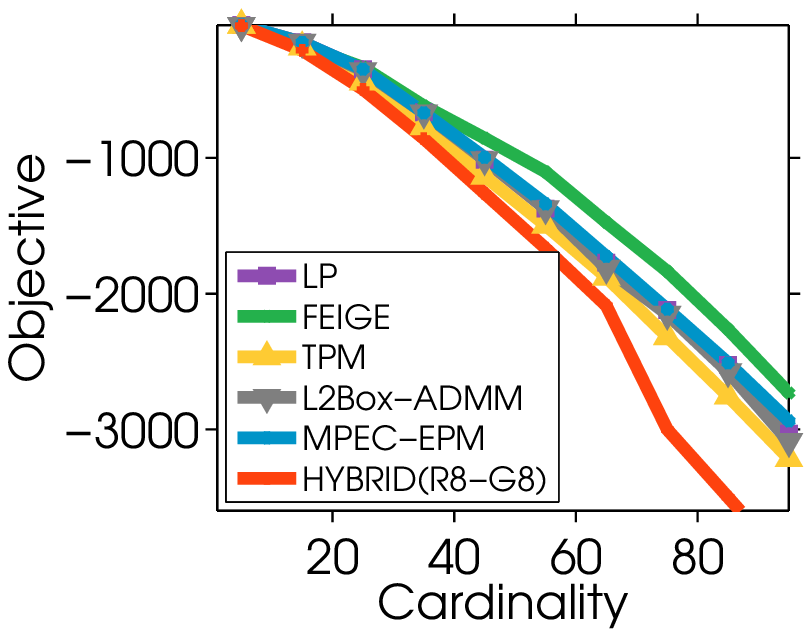}\vspace{-6pt} \caption{\footnotesize p2p-Gnutella08}\end{subfigure}\ghs
      \begin{subfigure}{\fourfigwid}\includegraphics[width=\objimgwid,height=\objimghei]{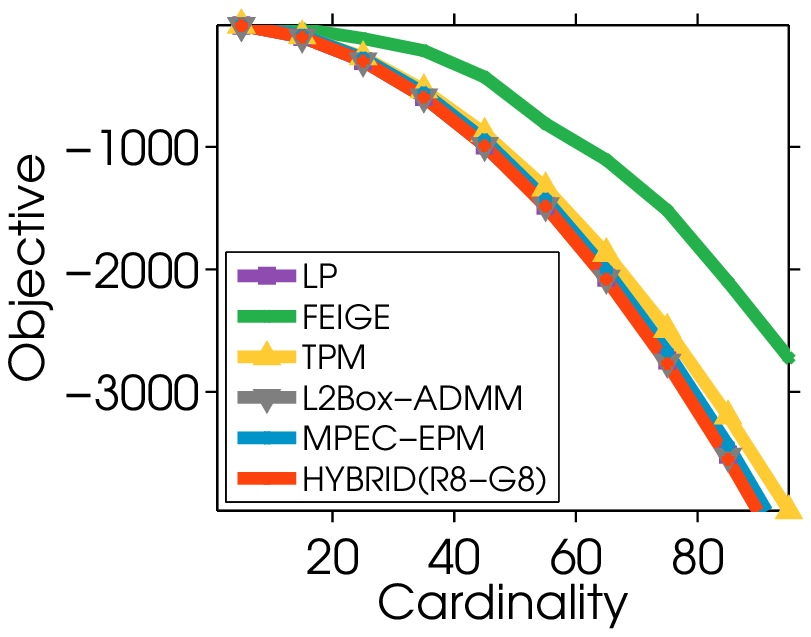}\vspace{-6pt} \caption{\footnotesize p2p-Gnutella09}\end{subfigure}\ghs
      \begin{subfigure}{\fourfigwid}\includegraphics[width=\objimgwid,height=\objimghei]{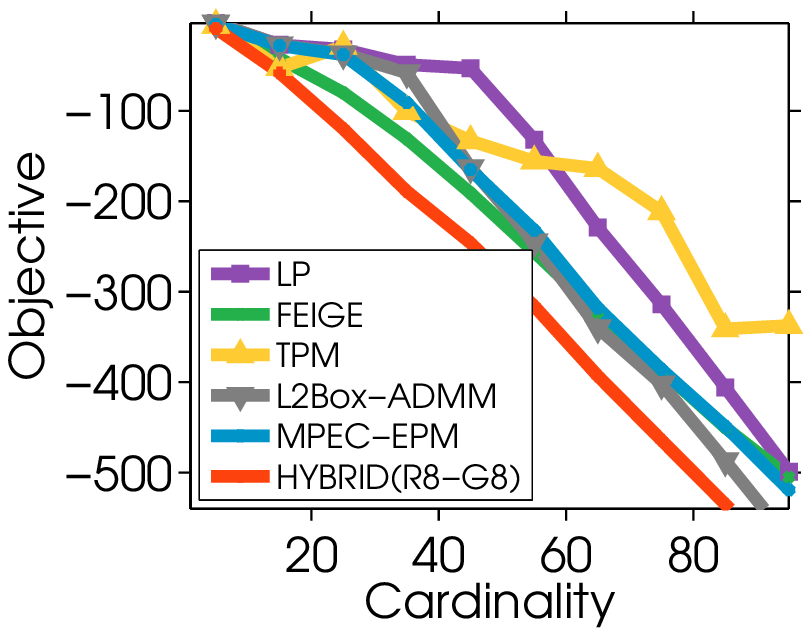}\vspace{-6pt} \caption{\footnotesize UK-2007-05}\end{subfigure}\\
\centering
\caption{Experimental results on dense subgraph discovery.}

\label{fig:subgraph}
\end{figure*}

\section{Dense Subgraph Discovery}
\label{sect:exp:dense}

Dense subgraph discovery is an important application of binary optimization. It aims at finding the maximum density subgraph on $s$ vertices \cite{ravi1994heuristic,feige2001dense,yuan2013truncated}, which can be formulated as the following binary program:
\beq \label{eq:subgraph}
\textstyle \min_{\bbb{x} \in \{0,1\}^n}~-\bbb{x}^T\bbb{W}\bbb{x},~s.t.~\bbb{x}^T\bbb{1} = s,
\eeq
\noi where $\bbb{W}\in \mathbb{R}^{n\times n}$ is the adjacency matrix of a graph. 


\bbb{Compared Methods.} We compare our HYBRID method on eight datasets (refer to the sub-captions in Figure \ref{fig:subgraph})\footnote{\url{https://snap.stanford.edu/data/}} against six methods: (i) LP relaxation \cite{YuanG17aaai}. (ii) Feige's greedy algorithm (GEIGE) \cite{feige2001dense}.  (iii) Truncated Power Method (TPM) \footnote{\url{https://sites.google.com/site/xtyuan1980/publications}} \cite{yuan2013truncated}. (iv) L2box-ADMM \cite{wu2016ellp}. (v) MPEC-EPM \cite{YuanG17aaai}. For more description of these methods, we refer to \cite{YuanG17aaai}. We show the average results of using 3 random initial points.

\bbb{Experimental Results.} Several observations can be drawn from Figure \ref{fig:subgraph}. (i) FEIGE generally fails to solve the dense subgraph discovery problem and it leads to solutions with low density. (ii) TPM gives better performance than state-of-the-art technique MPEC-EPM in some cases but it is unstable. (iii) Our proposed HYBRID generally outperforms all the compared methods.
\begin{figure*}{!ht}
\centering
\begin{subfigure}{\fourfigwid}\includegraphics[width=\objimgwid,height=\objimgheiB]{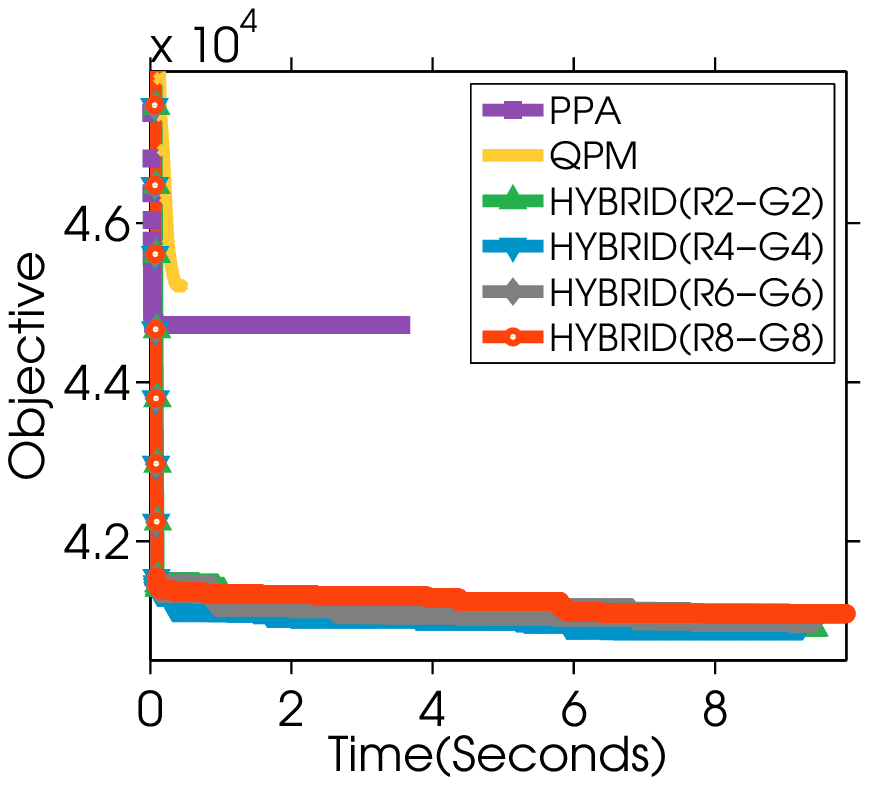}\vspace{-6pt} \end{subfigure}~~
\begin{subfigure}{\fourfigwid}\includegraphics[width=\objimgwid,height=\objimgheiB]{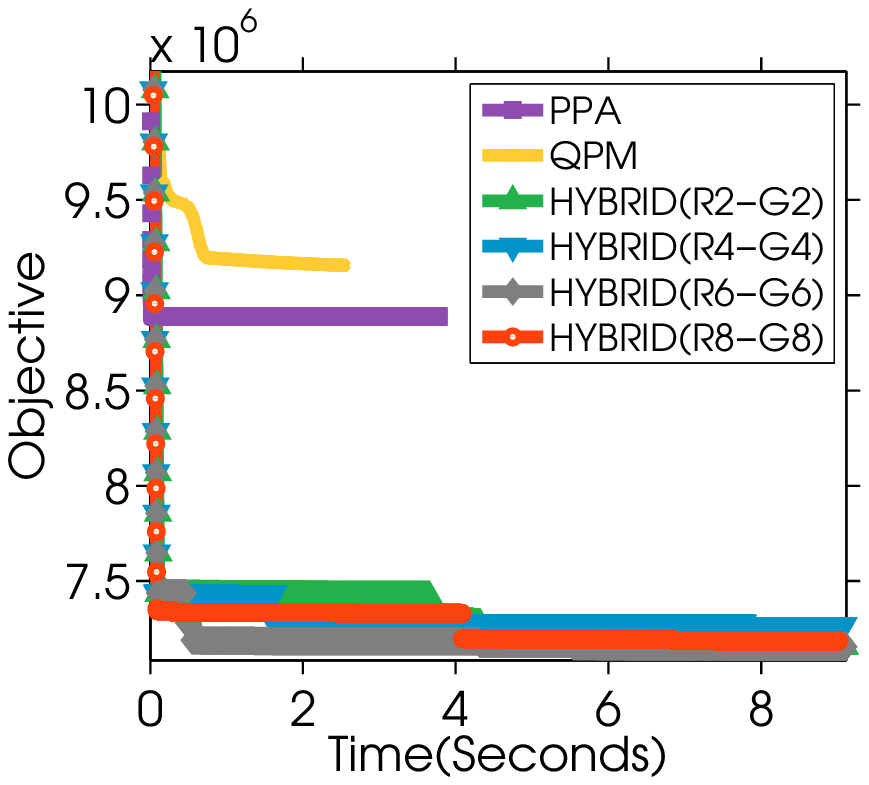}\vspace{-6pt} \end{subfigure}

\caption{Convergence curve of different methods for solving sparsity constrained least squares problems with $\{m,n,k\}=\{512, 2048, 20\}$. Left: `AI+bI'; Right: `AII+bII'.}
\label{fig:time}

\end{figure*}

\begin{figure*}[!t]
\centering
\begin{subfigure}{\fourfigwid}\includegraphics[width=\objimgwid,height=\objimghei]{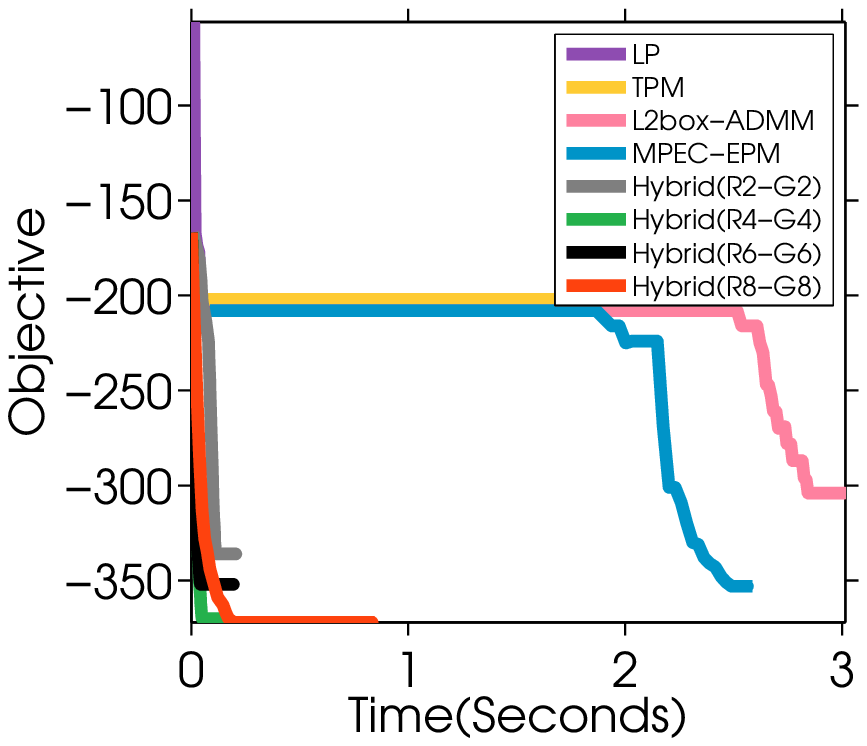}\vspace{-6pt} \end{subfigure}~~
\begin{subfigure}{\fourfigwid}\includegraphics[width=\objimgwid,height=\objimghei]{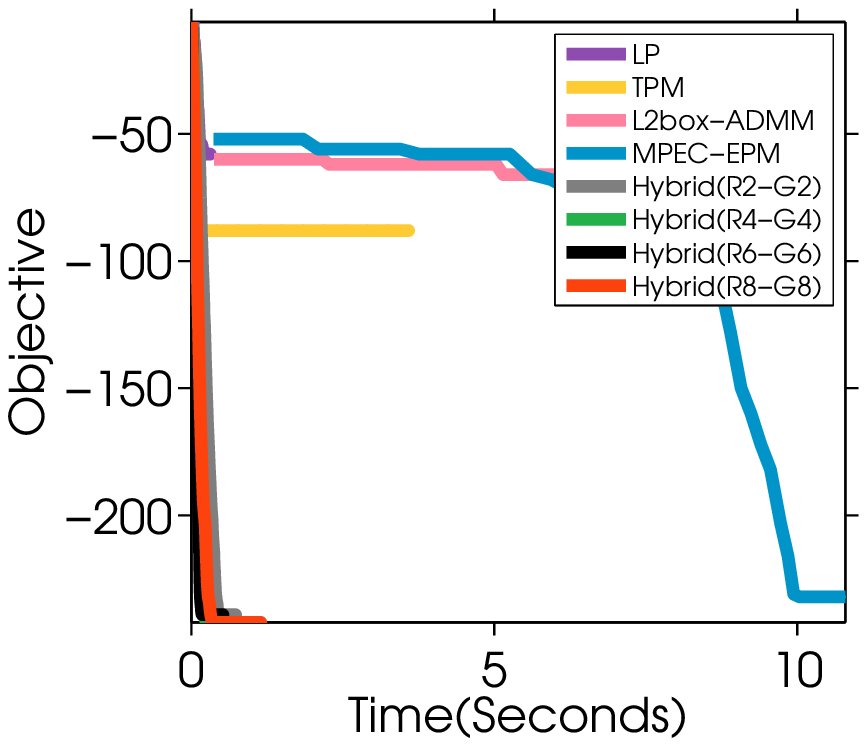}\vspace{-6pt} \end{subfigure}
\caption{Convergence curve of different methods for solving dense subgraph discovery problems. Left: `p2p-Gnutella09' data set; Right: `p2p-Gnutella30' data set. Theses two data sets contain 8114 and 36682 nodes, respectively. $k$ is set to 60 in our experiments.}
\label{fig:time:binary}

\end{figure*}

\subsection{Computational Efficiency of Algorithm \ref{algo:main}} \label{sect:computational:efficiency}

We show the convergence curve of different methods for binary optimization (see Figure \ref{fig:time}) and sparse optimization (see Figure \ref{fig:time:binary}). Generally speaking, our HYBRID is effective and practical for large-scale discrete optimization. Although it takes longer time to converge than PPA and QPM, the computational time is acceptable and it generally takes less than 30 seconds to converge in \emph{all} our instances. We think this computation time pays off as HYBRID achieves significantly higher accuracy than PPA and QPM. The main bottleneck of computation is on solving the small-sized subproblem using sub-exponential time $\mathcal{O}(2^k)$. The parameter $k$ in Algorithm \ref{algo:main} can be viewed as a tuning parameter to balance the efficacy and efficiency. One can further accelerate the algorithm using asynchronous parallelism or mini-batch optimization techniques.

\section{Conclusions}

This paper presents an effective and practical method for solving discrete optimization problems. Our method takes advantage of the effectiveness of combinatorial search and the efficiency of gradient descent. We also provided rigorous optimality analysis and convergence analysis for the proposed algorithm. The extensive experiments show that our method achieves state-of-the-art performance. 



\vspace{25pt}
{\huge Appendix}

\noi \appendix

The appendix section is organized as follows. Section \ref{sect:rate:binary}, \ref{sect:rate:sparse}, and \ref{sect:rate:sparse:c} present the convergence rate for binary optimization, sparse regularized optimization, and sparse constrained optimization, respectively.

\section{Convergence Rate for Binary Optimization} \label{sect:rate:binary}
We now prove the convergence rate of Algorithm \ref{algo:main} for binary optimization with $h \triangleq f_{\text{binary}}$. We define $ \Pi(\bbb{a}) \triangleq \arg\min_{\bbb{x}} \|\bbb{x}-\bbb{a}\|,~s.t.~ \bbb{x} \in \{-1,+1\}^n$. Our key technique is that when it holds that $\| \Pi(\bbb{x}) - \bbb{x}\|_2 \leq (1-\kappa)\|\Pi(\bar{\bbb{x}}) - \bbb{x}\|_2$ with $0<\kappa<1$ for any $\bbb{x}$, combining with the strongly convex property of $f(\cdot)$, we can establish global Q-linear convergence rate of Algorithm \ref{algo:main}. We remark that the strongly convexity assumption always holds for binary optimization, since one can append an additional term $\frac{\eta}{2}\|\bbb{x}\|_2^2$ to $f(\cdot)$ with sufficiently large $\eta$.

The following lemmas are useful in our proof.

\begin{lemma}
Assume that $f(\cdot)$ is $\alpha$-strongly convex. The following inequality holds for any $\bbb{x}$ and $\bbb{y}$:
\beq \label{eq:strongly:convex:function}
\frac{\alpha^2}{4} \|\bbb{y}-\bbb{x}\|_2^2 -\|\nabla f(\bbb{x})\|_2^2  \leq -\frac{\alpha}{2} (f(\bbb{x})-f(\bbb{y}))
\eeq
\begin{proof}

Since $f(\cdot)$ is $\alpha$-strongly convex, it holds that:
\beq \label{eq:alpha:strong}
f(\bbb{x}) - f(\bbb{y}) \leq \la \nabla f(\bbb{x}), \bbb{x} - \bbb{y}    \ra - \frac{\alpha}{2} \|\bbb{x} - \bbb{y}\|_2^2,~\forall \bbb{x},~\bbb{y}
\eeq

We naturally derive the following results:
\beq
\|\nabla f(\bbb{x})\|_2^2 - \frac{\alpha^2}{4} \|\bbb{y}-\bbb{x}\|_2^2&=&\left(\|\nabla f(\bbb{x})\|_2 - \frac{\alpha}{2} \|\bbb{y}-\bbb{x}\|_2\right) \cdot \left(\|\nabla f(\bbb{x})\|_2 + \frac{\alpha}{2} \|\bbb{y}-\bbb{x}\|_2\right)\nn\\
&\overset{}{\geq}&\left(\|\nabla f(\bbb{x})\|_2 - \frac{\alpha}{2} \|\bbb{y}-\bbb{x}\|_2\right) \cdot \left(0+\frac{\alpha}{2} \|\bbb{y}-\bbb{x}\|_2\right)\nn\\
&=& \frac{\alpha}{2} \left(\|\nabla f(\bbb{x})\|_2 \cdot \|\bbb{y}-\bbb{x}\|_2 - \frac{\alpha}{2} \|\bbb{y}-\bbb{x}\|_2^2\right)    \nn\\
&\overset{(a)}{\geq}&  \frac{\alpha}{2} \left( \la \nabla f(\bbb{x}), \bbb{y}-\bbb{x} \ra - \frac{\alpha}{2} \|\bbb{y}-\bbb{x}\|_2^2\right) \overset{(b)}{\geq} \frac{\alpha}{2} (f(\bbb{x}) - f(\bbb{y}) )\nn
\eeq
\noi where step $(a)$ uses the Cauchy-Schwarz inequality; step $(b)$ uses the $\alpha$-strongly convexity condition in (\ref{eq:alpha:strong}).

\end{proof}
\end{lemma}

\begin{lemma} \label{lemma:binary:local:error}
We define $\Pi(\bbb{a}) = \arg\min_{\bbb{x}} \|\bbb{x}-\bbb{a}\|,~s.t.~ \bbb{x} \in \{-1,+1\}^n$. The following inequality always holds for all $\bbb{x}$ and $\bbb{y}$:
\beq \label{eq:binary:local:error}
\| \Pi(\bbb{x}) - \bbb{x}\|_2^2 \leq (1-\kappa)\|\Pi(\bbb{y}) - \bbb{x}\|_2^2
\eeq
\noi with $\kappa=0$. In addition, if $\Pi(\bbb{x}_i) \neq \Pi(\bbb{y}_i)$ and $\bbb{x}_i \neq 0$ for some $i$, there exists a sufficiently small positive parameter $\kappa$ with $0<\kappa<1$ such that (\ref{eq:binary:local:error}) holds.

\begin{proof}
(i) First of all, noticing $\|\Pi(\bbb{x})\|_2^2=n$ for all $\bbb{x}$, we have the following results:
\beq
&&\| \Pi(\bbb{x}) - \bbb{x}\|_2^2 \leq (1-\kappa)\|\Pi(y) - \bbb{x}\|_2^2 \nn\\
&\Leftrightarrow&\kappa \|\Pi(y) - \bbb{x}\|_2^2 + \| \Pi(\bbb{x}) \|_2^2 + \|\bbb{x}\|_2^2 - 2 \la \Pi(\bbb{x}) , \bbb{x} \ra \leq \|\Pi(y)\|_2^2 + \|\bbb{x}\|_2^2 - 2\la \Pi(y),\bbb{x} \ra \nn\\
&\Leftrightarrow&\kappa \|\Pi(y) - \bbb{x}\|_2^2 \leq 2 \la \Pi(x) - \Pi(y),\bbb{x}\ra \nn
\eeq
\noi Since we have $\bbb{x}_i \cdot sign(\bbb{y}_i) \leq \bbb{x}_i \cdot sign(\bbb{x}_i)$ for all $i$, the conclusion of this lemma clearly holds when $\kappa=0$.

(ii) We now focus on the second part of this lemma. Note that $\Pi(\bbb{x}) \neq \Pi(\bbb{y})$ also implies $\bbb{x} \neq \Pi(\bbb{y})$ and we have $\|\Pi(y) - \bbb{x}\|_2^2>0$. Moreover,  Therefore, there exists a sufficient small $\kappa$ such that $\kappa \|\Pi(y) - \bbb{x}\|_2^2 \leq 2 \la \Pi(x) - \Pi(y),\bbb{x}\ra$ holds. This finishes the proof of this lemma.


\end{proof}
\end{lemma}

\begin{assumption} \label{assumption:binary}
Assuming that there exists a constant $0<\kappa<1$ such that:
\beq
\|\Pi(\bbb{x}) - \bbb{x}\|^2_2 \leq (1-\kappa) \|\bar{\bbb{x}} - \bbb{x}\|^2_2,~\forall \bbb{x}\nn
\eeq
\end{assumption}
\bbb{Remarks.} Based on the results in Lemma \ref{lemma:binary:local:error}, we conclude that when $\Pi(\bbb{x}_i) \neq \Pi(\bar{\bbb{x}}_i)$ and $\bbb{x}_i \neq 0$ for some $i$, such an assumption holds. This constant is similar to the classical local proximal error bound constant \cite{luo1993error} and it is necessary to use this constant to characterize the global linear convergence rate of the algorithm.

\setcounter{theorem}{1}
\begin{theorem}
\textbf{Proof of Convergence Rate when $h \triangleq f_{\text{binary}}$}. When Assumption 1 holds and $f(\cdot)$ is $s$-strongly convex, we have:
\beq \label{eq:binary:convergence}
F(\bbb{x}^{t}) - F(\bar{\bbb{x}}) \leq (1-\beta)^t (F(\bbb{x}^{0}) - F(\bar{\bbb{x}})),~~\text{with}~~\beta \triangleq(1-\kappa + \frac{\sqrt{{  (L+\theta-s)(1-\kappa) \kappa}}}{2(L+\theta)}) \frac{k}{n} < 1
\eeq
\noi In other words, it takes at most $\log_{(1-\beta)} (\frac{\epsilon}{F(\bbb{x}^0)-F(\bar{\bbb{x}})})$ times to find a solution $\bbb{x}^t$ satisfying $F(\bbb{x}^t) - F(\bar{\bbb{x}}) \leq \epsilon$.

\begin{proof}
First of all, we define
\beq
 \textstyle d(\bbb{x}) \triangleq \arg \min_{\bbb{z}\in \mathbb{R}^n} H(\bbb{x},\bbb{z}),~H(\bbb{x},\bbb{z}) \triangleq f(\bbb{x}) + \la \nabla f(\bbb{x}), \bbb{z} \ra + \frac{L+\theta}{2} \|\bbb{z}\|_2^2 + h_{\text{binary}}(\bbb{x}+\bbb{z}) \label{eq:ddd2}. \nn
\eeq
\noi Since $\bbb{x}^t \in \Psi$ for all $t$. We have the following inequalities:
\beq \label{eq:stoc:bound}
 \E [f(\bbb{x}^{t+1})] &\overset{(a)}{\leq}& \textstyle \E [f(\bbb{z})] + \frac{\theta}{2} \|\bbb{z}-\bbb{x}^t\|_2^2 - \frac{\theta}{2} \|\bbb{x}^{t+1} - \bbb{x}^t\|_2^2,~\forall \bbb{z}_N = \bbb{x}_N^t,~\bbb{z} \in  \Psi \nn\\
&\overset{(b)}{\leq}& \textstyle \E [f(\bbb{z})] + \frac{\theta}{2} \|\bbb{z}-\bbb{x}^t\|_2^2 ,~\forall \bbb{z}_N = \bbb{x}_N^t,~\bbb{z} \in  \Psi \nn\\
&\overset{(c)}{\leq}& \textstyle \E [f(\bbb{x}^t) + \la \nabla f(\bbb{x}^t),~ \bbb{z} - \bbb{x}^t \ra + \frac{L}{2}\| \bbb{z} - \bbb{x}^t \|_2^2 + \frac{\theta}{2} \|\bbb{z}-\bbb{x}^t\|_2^2] ,~\forall \bbb{z}_N = \bbb{x}_N^t,~\bbb{z} \in  \Psi \nn\\
&\overset{(d)}{=}& \textstyle f(\bbb{x}^t) + \E [\la (\nabla f(\bbb{x}^t))_B,~ \bbb{z}_B - \bbb{x}_B^t \ra] + \frac{L+\theta}{2}\E [\| \bbb{z}_B - \bbb{x}^t_B \|_2^2 ] ,~\bbb{z} \in  \Psi \nn\\
&\overset{(e)}{\leq}& \textstyle f(\bbb{x}^t) + \E [\la (\nabla f(\bbb{x}^t))_B,~ (\bbb{x}^t+d(\bbb{x}^t))_B - \bbb{x}_B^t \ra] + \frac{L+\theta}{2}\E [\| (\bbb{x}^t+d(\bbb{x}^t))_B - \bbb{x}^t_B \|_2^2 ] ,~\bbb{z} \in  \Psi \nn\\
&\overset{(f)}{=}& \textstyle f(\bbb{x}^t) + \frac{k}{n} \la \nabla f(\bbb{x}^t),~ d(\bbb{x}^t)\ra + \frac{k}{n}\frac{L+\theta}{2} \| d(\bbb{x}^t)\|_2^2
\eeq
\noi where step $(a)$ uses the optimality condition of $\bbb{x}^{t+1}$; step $(b)$ uses the fact that $- \frac{\theta}{2} \|\bbb{x}^{t+1} - \bbb{x}^t\|_2^2\leq 0$; step $(c)$ uses the Lipschitz continuity of the gradient of $f(\cdot)$ that: $f(\bbb{z}) \leq f(\bbb{x}) + \la \nabla f(\bbb{x}),~\bbb{z}-\bbb{x} \ra + \frac{L}{2}\|\bbb{x}-\bbb{z}\|_2^2,~\forall \bbb{x},~\bbb{z}$; step $(d)$ uses the fact that $\bbb{z}_N = \bbb{x}_N^t$; step (e) uses the choice that $\bbb{z} = \bbb{x} + d(\bbb{x})$; step $(f)$ uses the fact that each block $B$ is picked randomly with probability $k/n$.

We denotes
\beq \label{def:x:breve}
\breve{\bbb{x}}^t \triangleq \bbb{x}^t - \nabla f(\bbb{x}^t)/(L+\theta).
\eeq
\noi For any $\bbb{x}^{t} \in \Psi$ and $\bbb{x}^{t+1} \in \Psi$, we naturally derive the following inequalities:
\beq \label{eq:binary:one}
&& \E[\frac{n}{k} F(\bbb{x}^{t+1})- \frac{n}{k} F(\bbb{x}^t)] \nn\\
&\overset{(a)}{\leq} & \E[ \la \nabla f(\bbb{x}^t), d(\bbb{x}^t) \ra + \frac{L+\theta}{2} \|d(\bbb{x}^t)\|_2^2] \nn\\
&\overset{(b)}{=} &  \E[\la \nabla f(\bbb{x}^t), \Pi\left( \bbb{x}^t - \nabla f(\bbb{x}^t)/(L+\theta)\right)-\bbb{x}^t \ra + \frac{L+\theta}{2} \|\Pi\left(\bbb{x}^t - \nabla f(\bbb{x}^t)/(L+\theta)\right)-\bbb{x}^t\|_2^2] \nn\\
&\overset{(c)}{=} &  \E[\la \nabla f(\bbb{x}^t), \Pi ( \breve{\bbb{x}}^t  )-\bbb{x}^t \ra + \frac{L+\theta}{2} \|\Pi ( \breve{\bbb{x}}^t  )-\bbb{x}^t\|_2^2] \nn\\
&\overset{(d)}{=}& \E[\frac{L+\theta}{2} \|\Pi\left( \breve{\bbb{x}}^t \right)- \breve{\bbb{x}}^t \|_2^2  -  \frac{L+\theta}{2}\|\nabla f(\bbb{x}^t)/(L+\theta)\|_2^2 ]
\eeq
\noi where step $(a)$ uses (\ref{eq:stoc:bound}); step $(b)$ uses the fact that $d(\bbb{x})$ has the following closed-form solution that $d(\bbb{x}) = \Pi\left(\bbb{x} - \frac{\nabla f(\bbb{x})}{L+\theta}\right)-\bbb{x}$ for any $\bbb{x}$; step $(c)$ uses the definition of $\breve{\bbb{x}}^t$ in (\ref{def:x:breve}); step $(d)$ the inequality that $\frac{t}{2}\|\bbb{x}\|_2^2 +\la \bbb{x},~\bbb{a}\ra = \frac{t}{2}\|\bbb{x}+\bbb{a}/t\|_2^2 - \frac{t}{2}\|\bbb{a}/t\|_2^2 $ for all $\bbb{a},~\bbb{x},~t$.

Based on Assumption \ref{assumption:binary}, we have the following important inequality:
\beq\label{eq:important:assumption}
\|\Pi\left( \breve{\bbb{x}}^t \right)- \breve{\bbb{x}}^t \|_2^2 \leq (1-\kappa)\|\Pi(\bar{\bbb{x}}) - \breve{\bbb{x}}^t \|_2^2
\eeq

We now bound the inequality in (\ref{eq:binary:one}). We derive the following results:
\beq\label{eq:binary:two}
&& \E[\frac{n}{k} F(\bbb{x}^{t+1})- \frac{n}{k} F(\bbb{x}^t)] \nn\\
&\overset{(a)}{\leq}& \E[\frac{  (L+\theta)(1-\kappa)}{2}  \|\Pi(\bar{\bbb{x}}) - \breve{\bbb{x}}^t \|_2^2 -  \frac{L+\theta}{2}\|\nabla f(\bbb{x}^t)/(L+\theta)\|_2^2 ] \nn\\
&\overset{(b)}{=}& \E[\frac{  (L+\theta)(1-\kappa)}{2}  \|\bar{\bbb{x}} -  \bbb{x}^t + \nabla f(\bbb{x}^t)/(L+\theta) \|_2^2 -  \frac{L+\theta}{2}\|\nabla f(\bbb{x}^t)/(L+\theta)\|_2^2 ] \nn\\
&\overset{}{=}& \E[\frac{  (L+\theta)(1-\kappa)}{2}  \|\bar{\bbb{x}} - \bbb{x}^t \|_2^2 + (1-\kappa) \la \bar{\bbb{x}} - \bbb{x}^t, \nabla f(\bbb{x}^t)  \ra  -\frac{\kappa}{2} \|\nabla f(\bbb{x}^t)/(L+\theta)\|_2^2] \nn\\
&\overset{(c)}{\leq}& \E[\frac{  (L+\theta-s)(1-\kappa)}{2}  \|\bar{\bbb{x}} - \bbb{x}^t \|_2^2 + (1-\kappa) ( f(\bar{\bbb{x}}) - f(\bbb{x}^t)   )  -\frac{\kappa}{2} \|\nabla f(\bbb{x}^t)/(L+\theta)\|_2^2] \nn\\
&\overset{(d)}{=}& \E[(\frac{\kappa}{2(L+\theta)^2}) (\frac{  (L+\theta-s)(1-\kappa)(L+\theta)^2}{ \kappa}  \|\bar{\bbb{x}} - \bbb{x}^t \|_2^2 - \|\nabla f(\bbb{x}^t)\|_2^2) + (1-\kappa) ( f(\bar{\bbb{x}}) - f(\bbb{x}^t)  ) ] \nn\\
&\overset{(e)}{\leq} & \E[ (\frac{\kappa}{2(L+\theta)^2}) \sqrt{\frac{  (L+\theta-s)(1-\kappa)(L+\theta)^2}{ \kappa}  } (- f(\bbb{x}^t)  + f(\bar{\bbb{x}})  ) + (1-\kappa) ( f(\bar{\bbb{x}}) - f(\bbb{x}^t)  ) ] \nn \\
&\overset{}{=} & \E[ (1-\kappa + \frac{\sqrt{{  (L+\theta-s)(1-\kappa) \kappa}}}{2(L+\theta)})( f(\bar{\bbb{x}}) - f(\bbb{x}^t)  ) ]
\eeq
\noi step $(a)$ uses (\ref{eq:important:assumption}); step $(b)$ uses the fact that $\Pi(\bar{\bbb{x}})=\bar{\bbb{x}}$ and the definition of $\breve{\bbb{x}}^t$ in (\ref{def:x:breve}); step $(c)$ uses the $s$-strongly convexity of $f(\cdot)$ that $\la \bar{\bbb{x}}-\bbb{x}^t, \nabla f(\bbb{x}^t)\ra \leq f(\bar{\bbb{x}}) - f(\bbb{x}^t) - \frac{s}{2}\|\bar{\bbb{x}}-\bbb{x}^t\|_2^2$; step $(e)$ uses the inequality in (\ref{eq:strongly:convex:function}). It it worthwhile to note that when Assumption \ref{assumption:binary} does not hold with $\kappa\rightarrow 0$, step $(d)$ breaks down.

Based on (\ref{eq:binary:two}) and the definition of $\beta$ in (\ref{eq:binary:convergence}), we have the following inequality: $\E[F(\bbb{x}^{t+1}) - F(\bbb{x}^t)] \leq \beta \E[ F(\bar{\bbb{x}}) - F(\bbb{x}^t)]$. Rearranging terms, we obtain that $\E[F(\bbb{x}^{t+1}) - F(\bar{\bbb{x}})] \leq \E[(1-\beta) (F(\bbb{x}^t)-F(\bar{\bbb{x})})]$. In other words, the sequence $\{f(\bbb{x}^t)\}$ converges to the stationary point linearly in the quotient sense. Solving this recursive formulation, we obtain (\ref{eq:binary:convergence}). Therefore, we conclude that it takes at most $\log_{(1-\beta)} (\frac{\epsilon}{F(\bbb{x}^0)-F(\bar{\bbb{x}})})$ times  to find a local optimal solution satisfying $F(\bbb{x}^t) - F(\bar{\bbb{x}}) \leq \epsilon$.

\end{proof}

\end{theorem}

\section{Convergence Rate for Sparse Regularized Optimization}
\label{sect:rate:sparse}

We now prove the convergence rate of Algorithm \ref{algo:main} for sparse regularized optimization with $h \triangleq f_{\text{sparse}}$. We have derived the upper bound for the number of changes $J$ for the support set in Theorem \ref{theorem:convergence}. We now need to derive a bound on the number of iterations performed after the support set is fixed. By combing these two bounds, we establish our proof.

The following lemma is useful in our proof.

\begin{lemma} \label{eq:recursive}
 Assume a nonnegative sequence $\{u^t\}_{t=0}^{\infty}$ satisfies $(u^{t+1})^2 \leq C (u^{t} - u^{t+1})$ for some constant $C$. We have:
\beq
u^{t} \leq \frac{\max(2C, \sqrt{Cu^0})}{t}
\eeq

\begin{proof}

We denote $C_1 \triangleq \max(2C, \sqrt{Cu^0})$. Solving this quadratic inequality, we have:
\beq
u^{t+1} \leq -  \frac{C}{2} +  \frac{C}{2} \sqrt{1+\frac{ 4 u^t}{C}}
\eeq
\noi We now show that $u^{t+1} \leq  \frac{C_1}{t+1}$, which can be obtained by mathematical induction. (i) When $t=0$, we have $u^{1}  \leq - \frac{C}{2} + \frac{C}{2} \sqrt{1+ \frac{4 u^0}{C} } \leq -\frac{C}{2} + \frac{C}{2} (1+\sqrt{\frac{4u^0}{C}} ) = \frac{C}{2}\sqrt{\frac{4u^0}{C}} = \sqrt{Cu^0 } \leq \frac{C_1}{t+1}$. (ii) When $t\geq 1$, we assume that $u^{t} \leq  \frac{C_1}{t}$ holds. We derive the following results: $t\geq 1 \Rightarrow \frac{t+1}{t} \leq 2$ $~\overset{(a)}{\Rightarrow}~  C \frac{t+1}{t} \leq C_1$ $~\overset{(b)}{\Rightarrow}~ $ $ C (\frac{1}{t} - \frac{1}{t+1} ) \leq \frac{C_1}{(t+1)^2}$ $\Rightarrow \frac{C}{t} \leq    \frac{C}{t+1} + \frac{C_1}{(t+1)^2}$$\Rightarrow \frac{C C_1}{t} \leq    \frac{C C_1 }{t+1} + \frac{C^2_1}{(t+1)^2}$$\Rightarrow \frac{C^2}{4} + \frac{C C_1}{t} \leq    \frac{C C_1 }{t+1} + \frac{C^2_1}{(t+1)^2}+\frac{C^2}{4}$$\Rightarrow \frac{C^2}{4} (1+ \frac{4 C_1 }{Ct}  )  \leq   (\frac{C}{2} + \frac{C_1}{t+1})^2$$\Rightarrow \frac{C}{2} \sqrt{1+ \frac{4 C_1 }{Ct} }  \leq   \frac{C}{2} + \frac{C_1}{t+1}$$\Rightarrow - \frac{C}{2} +  \frac{C}{2} \sqrt{1+ \frac{4 C_1 }{Ct} }  \leq    \frac{C_1}{t+1}$$~\overset{(c)}{\Rightarrow}~ - \frac{C}{2} +  \frac{C}{2} \sqrt{1+ \frac{4 u^t }{C} }  \leq    \frac{C_1}{t+1}$ $\Rightarrow   u^{t+1} \leq \frac{C_1}{t+1}$. Here, step $(a)$ uses $2C \leq C_1$; step $(b)$ uses $\frac{1 }{t(t+1)} =  \frac{1}{t} - \frac{1}{t+1} $; step $(c)$ uses $u^{t} \leq  \frac{C_1}{t}$.

\end{proof}

\end{lemma}

The following proposition establishes a bound on the number of iterations performed after the support set is fixed, which is novel in this paper. Note that the optimization problems become convex with fixing support set, and any stationary point $\bar{\bbb{x}}$ is also the global optimal solution for the convex problems.

\begin{proposition} \label{the:sparse:bound:strong:yn}
Assume that the support set of $\bbb{x}^t$ does not changes for all $t$. We have the following results:

(i) When $f(\cdot)$ is convex, it takes at most $V_1$ iterations in expectation for Algorithm \ref{algo:main} to converge to a stationary point $\bar{\bbb{x}}$ satisfying $F(\bbb{x}^t)-F(\bar{\bbb{x}}) \leq \epsilon$, where $V_1$ is defined as:
\beq \label{eq:convex:boundt}
V_1 = \frac{\max( \frac{4 \nu^2}{\theta} , \sqrt{  \frac{2 \nu^2 [F(\bbb{x}^0) - F(\bar{\bbb{x}})]}{\theta } })}{\epsilon },~\text{with}~\nu\triangleq \frac{2n \rho \sqrt{k}\theta}{k}.
\eeq

(ii) When $f(\cdot)$ is $s$-strongly convex, it takes at most $V_2$ iterations in expectation for Algorithm \ref{algo:main} to converge to a stationary point $\bar{\bbb{x}}$ satisfying $F(\bbb{x}^t)-F(\bar{\bbb{x}}) \leq \epsilon$, where $V_2$ is defined as:
\beq \label{eq:convex:boundt2}
V_2 = \log_{\frac{\varpi}{1+\varpi}} (\frac{\epsilon}{F(\bbb{x}^{0}) - F(\bar{\bbb{x}})}),~\text{with}~\varpi\triangleq \frac{\theta n}{sk}.
\eeq

\begin{proof}

First of all, we notice that when the support set is fixed, the original problem reduces to the following convex composite optimization problem:
\beq
\textstyle F(\bbb{x}) \triangleq f(\bbb{x}) + p(\bbb{x}) + const,~~~\text{with}~~~ p(\bbb{x}) \triangleq I_{\Omega}(\bbb{x}) \label{eq:reduced:convex}
\eeq

\noi We use $\partial F(\bbb{x})$ to denote the sub-gradient of $F(\cdot)$ in $\bbb{x}$. Since the algorithm solves the following optimization: $\bbb{x}^{t+1} \Leftarrow \arg \min_{\bbb{z}}~f(\bbb{z}) + p(\bbb{z}) + \frac{\theta}{2} \|\bbb{z}-\bbb{x}^{t}\|^2,~s.t.~\bbb{z}_{N} = \bbb{x}^t_{N}$, we have the following optimality condition for $\bbb{x}^{t+1}$:
\beq \label{eq:optimality:sparse}
\bbb{w}_{B} + \theta (\bbb{x}^{t+1}-\bbb{x}^{t})_B = \bbb{0},~\bbb{x}^{t+1}_N = \bbb{x}^t_N,~\forall \bbb{w} \in \partial F(\bbb{x}^{t+1})
\eeq

(i) We now consider the case when $f(\cdot)$ is generally convex. We derive the following inequalities:
\beq
\E[F(\bbb{x}^{t+1})] - F(\bar{\bbb{x}})&\overset{(a)}{\leq}& \E[\la\bbb{w} ,  \bbb{x}^{t+1} - \bar{\bbb{x}} \ra],~\forall \bbb{w} \in  \partial F(\bbb{x}^{t+1}) \nn\\
&\overset{(b)}{\leq}& \E[\frac{n}{k} \la \bbb{w}_{B} , (\bbb{x}^{t+1} - \bar{\bbb{x}})_{B} \ra],~\forall \bbb{w} \in \partial F(\bbb{x}^{t+1}) \nn\\
&\overset{(c)}{=}& \E[\frac{n}{k} \la - \theta (\bbb{x}^{t+1}-\bbb{x}^{t})_B , (\bbb{x}^{t+1} - \bar{\bbb{x}})_{B} \ra]\nn\\
&\overset{(d)}{\leq} & \E[\frac{n}{k}   \theta \|(\bbb{x}^{t+1}-\bbb{x}^{t})_B\|_2 \cdot \|(\bbb{x}^{t+1} - \bar{\bbb{x}})_{B}\|_2] \nn\\
&\overset{(e)}{=} & \E[\frac{n}{k}   \theta \|(\bbb{x}^{t+1}-\bbb{x}^{t})\|_2 \cdot \|(\bbb{x}^{t+1} - \bar{\bbb{x}})_{B}\|_2] \nn\\
&\overset{(f)}{\leq} &  \E[\underbrace{\frac{n}{k}   \theta \cdot {2\rho}{\sqrt{k}}}_{\nu} \|\bbb{x}^{t+1}-\bbb{x}^{t}\|_2]  \label{eq:cost:to:go}
\eeq
\noi where step $(a)$ uses the convexity of $F$; step $(b)$ uses the fact that each block $B$ is picked randomly with probability $k/n$; step $(c)$ uses the optimality condition of $\bbb{x}^{t+1}$ in (\ref{eq:optimality:sparse}); step $(d)$ uses the Cauchy-Schwarz inequality; step $(e)$ uses $\|(\bbb{x}^{t+1}-\bbb{x}^{t})\|_2 = \|(\bbb{x}^{t+1}-\bbb{x}^{t})_B\|_2$; step $(f)$ uses $\|(\bbb{x}^{t+1} - \bar{\bbb{x}})_B\| \leq {\sqrt{k}}\|(\bbb{x}^{t+1} - \bar{\bbb{x}})_B\|_{\infty}\leq {\sqrt{k}}\|\bbb{x}^{t+1} - \bar{\bbb{x}}\|_{\infty}\leq {\sqrt{k}}(\|\bbb{x}^{t+1}\|_{\infty} + \|\bar{\bbb{x}}\|_{\infty})\leq {2 \rho}{\sqrt{k}}$.

For any $\bbb{x}^t,\bbb{x}^{t+1},~\bar{\bbb{x}}~\in\Omega$, we derive the following results:
\beq
\E[F(\bbb{x}^{t+1}) - F(\bar{\bbb{x}})] \overset{(a)}{\leq} \E[\nu \|\bbb{x}^{t+1}-\bbb{x}^t\|_2] \overset{(b)}{\leq} \E[\nu \sqrt{ \frac{2}{\theta} \left(F(\bbb{x}^{t}) - F(\bbb{x}^{t+1})\right)    }]
\eeq
\noi where the step $(a)$ uses (\ref{eq:cost:to:go}); step $(b)$ uses the sufficient decent condition in (\ref{eq:suff:dec}). Denoting $\Delta^{t} \triangleq \E[F(\bbb{x}^{t}) - F(\bar{\bbb{x}})]$ and $C \triangleq \frac{2 \nu^2}{\theta}$, we have the following inequality:
\beq
(\Delta^{t+1})^2 \leq C (\Delta^{t} - \Delta^{t+1})\nn
\eeq
\noi Combining with Lemma \ref{eq:recursive}, we have:
\beq
\E[F(\bbb{x}^{t}) - F(\bar{\bbb{x}})]  \leq \frac{\max( \frac{4 \nu^2}{\theta} , \sqrt{  \frac{2 \nu^2 \Delta^0}{\theta } })}{t} \nn
\eeq
\noi Therefore, we obtain (\ref{eq:convex:boundt}).

(ii) We now consider the case when $f(\cdot)$ is strongly convex. For any $\bbb{x}^t,\bbb{x}^{t+1},~\bar{\bbb{x}}~\in\Omega$, we derive the following results:
\beq \label{eq:stronlgy:convex:J2}
\E[F(\bbb{x}^{t+1}) - F(\bar{\bbb{x}})   ] &\overset{(a)}{\leq}& \E[-\frac{s}{2} \|\bar{\bbb{x}}-\bbb{x}^{t+1}\|_2^2 - \la \bar{\bbb{x}}-\bbb{x}^{t+1},  \bbb{w} \ra],~\forall \bbb{w} \in \partial F(\bbb{x}^{t+1}) \nn\\
&\overset{(b)}{\leq}& \E[\frac{1}{2s} \|\bbb{w} \|_2^2],~\forall \bbb{w} \in \partial F(\bbb{x}^{t+1}) \nn\\
&\overset{(c)}{=}&\E[ -\frac{1}{2s}  \|\bbb{w}_{B}  \|_2^2 \times \frac{n}{k}],~\forall \bbb{w} \in \partial F(\bbb{x}^{t+1}) \nn\\
&\overset{(d)}{=}& \E[-\frac{1}{2s}  \|   \theta (\bbb{x}^{t+1}-\bbb{x}^{t})_B   \|_2^2 \times  \frac{n}{k}] \nn\\
&\overset{}{=}& \E[  \frac{\theta^2 n}{2sk} \| \bbb{x}^{t+1}-\bbb{x}^{t}\|_2^2]  \nn\\
&\overset{(e)}{\leq}& \E[   \frac{\theta^2 n}{2sk} \frac{2}{\theta} \left(F(\bbb{x}^{t}) - F(\bbb{x}^{t+1})]\right) \nn\\
& \overset{}{=}&\E[  \underbrace{\frac{\theta n}{sk}}_{\varpi} \left( [F(\bbb{x}^{t}) - F(\bar{\bbb{x}}) ] -  [F(\bbb{x}^{t+1}) - F(\bar{\bbb{x}})]   \right) ]
\eeq
\noi where step $(a)$ uses the strongly convexity of $f(\cdot)$; step $(b)$ uses the fact that
$-\frac{s}{2}\|\bbb{x}\|_2^2 - \la \bbb{x},\bbb{y}\ra \leq \frac{1}{2s}\|\bbb{y}\|_2^2$; step $(c)$ uses the fact that $\E[\| \bbb{w}_{B}\|_2^2] = \frac{k}{n}\|\bbb{w}\|_2^2$; step $(d)$ uses the optimality of $\bbb{x}^{t+1}$; step $(e)$ uses the sufficient condition in $(\ref{eq:suff:dec})$.

Rearranging terms for (\ref{eq:stronlgy:convex:J2}), we have: $\frac{\E[F(\bbb{x}^{t+1}) - F(\bar{\bbb{x}})]}{\E[F(\bbb{x}^{t}) - F(\bar{\bbb{x}})]} \leq \frac{\varpi}{1+\varpi}$. Solving the recursive formulation, we obtain:
\beq \label{eq:aaa}
\E[F(\bbb{x}^{t}) - F(\bar{\bbb{x}})] \leq \E[\left(\frac{\varpi}{1+\varpi}\right)^t [F(\bbb{x}^{0}) - F(\bar{\bbb{x}})]], \nn
\eeq
\noi and it holds that $t \leq \log_{\frac{\varpi}{1+\varpi}} \left(\frac{F(\bbb{x}^{t}) - F(\bar{\bbb{x}})}{F(\bbb{x}^{0}) - F(\bar{\bbb{x}})}\right)$ in expectation. Therefore, we obtain (\ref{eq:convex:boundt2}).

\end{proof}
\end{proposition}

\setcounter{theorem}{2}

\begin{theorem}
\textbf{Convergence Rate when $h \triangleq f_{\text{sparse}}$}. We have the follow results:

(i) When $f(\cdot)$ is generally convex, it takes at most $N_1$ iterations in expectation to converge to a block-$k$ stationary point $\bar{\bbb{x}}$ such that $F(\bbb{x}^t)-F(\bar{\bbb{x}})\leq \epsilon$, where
\beq
N_1=\left( \frac{\bar{J}}{D} + \frac{1}{\epsilon}\right) \times  \max\left( \frac{4 \nu^2}{\theta} , \sqrt{  \frac{2 \nu^2 (F(\bbb{x}^0) - F(\bar{\bbb{x}}) - D)}{\theta } }\right) .\nn
\eeq

(ii) When $f(\cdot)$ is $s$-strongly convex, it takes at most $N_2$ iterations in expectation to converge to a block-$k$ stationary point $\bar{\bbb{x}}$ such that $F(\bbb{x}^t)-F(\bar{\bbb{x}})\leq \epsilon$, where
\beq
N_2=  \bar{J} \log_{\frac{\varpi}{1+\varpi}} \left( \frac{D}{(F(\bbb{x}^0)-F(\bar{\bbb{x}}))} \right)  + \log_{\frac{\varpi}{1+\varpi}} \left(\frac{\epsilon}{F(\bbb{x}^{0}) - D - F(\bar{\bbb{x}})}\right).
\eeq

\begin{proof}

(i) We first consider the case when $f(\cdot)$ is generally convex. We denote $Z_t = \{i:\bbb{x}^t_i=0\}$. We known that the $Z_t$ only changes for a finite number of times. We assume that $Z_t$ only changes at $t=c_1,c_2,...,c_{\bar{J}}$ and we define $c_0=0$. Therefore, we have:
\beq
Z_{0} = Z_{1}=,...,Z_{-1+c_{1}} \neq Z_{c_{1}}=Z_{1+c_{1}}=Z_{2+c_{1}} = ,...,=Z_{-1+c_{j}} \neq Z_{c_{j}} = ... \neq  Z_{c_{{\bar{J}}}} = ...\nn
\eeq
\noi with $j=1,...,{\bar{J}}$. We denote $\bar{\bbb{x}}^{c_j}$ as the optimal solution of the following optimization problem:
\beq \label{eq:convex:subprob}
\min_{\bbb{x}}~f(\bbb{x}) + p(\bbb{x}),~s.t.~\bbb{x}_{I_{c_j}} = 0
\eeq
\noi with $1\leq j\leq {\bar{J}}$.

The solution $\bbb{x}^{c_j}$ changes $j$ times, the objective values decrease at least by $j D$, where $D$ is defined in (\ref{eq:DD}). Therefore, we have:
\beq
F(\bbb{x}^{c_j}) \leq F(\bbb{x}^0) - {j}\times D \nn
\eeq
\noi Combing with the fact that $ F(\bar{\bbb{x}}) \leq F(\bar{\bbb{x}}^{c_j})$, we obtain:
\beq \label{eq:J:bound00}
0\leq F(\bbb{x}^{c_{j}}) - F(\bar{\bbb{x}}^{c_j}) \leq F(\bbb{x}^0) - F(\bar{\bbb{x}}) - j \times D
\eeq
\noi We now focus on the intermediate solutions $\bbb{x}_{c_{j-1}},~\bbb{x}_{1+c_{j-1}},...,~\bbb{x}_{-1+c_j},\bbb{x}_{c_j}$. Using part (i) in Proposition \ref{the:sparse:bound:strong:yn}, we conclude that to obtain an accuracy such that $ F(\bbb{x}^{c_{j}}) - F(\bar{\bbb{x}}^{c_j}) \leq  D$, it takes at most ${ \max\left( \frac{4 \nu^2}{\theta} , \sqrt{  \frac{2 \nu^2 (F(\bbb{x}^{c_{j}}) - F(\bar{\bbb{x}}^{c_j}))}{\theta } }\right)   }/{D}$ iterations to converge to $\bar{\bbb{x}}^{c_j}$, that is,

\beq \label{eq:J:bound1}
c_j-c_{j-1} \leq \frac{\max\left( \frac{4 \nu^2}{\theta} , \sqrt{  \frac{2 \nu^2 (F(\bbb{x}^{c_{j}}) - F(\bar{\bbb{x}}^{c_j}))}{\theta } }\right)}{D} \leq \frac{\max\left( \frac{4 \nu^2}{\theta} , \sqrt{  \frac{2 \nu^2 (F(\bbb{x}^0) - F(\bar{\bbb{x}}) - j \times D)}{\theta } }\right)}{D}
 \eeq


\noi Summing up the inequality in (\ref{eq:J:bound1}) for $j = 1,2,...,\bar{J}$ and using the fact that $j\geq 1$ and $c_0=0$, we obtain that:
\beq \label{eq:J:total}
c_{{\bar{J}}} \leq  \bar{J} \times \frac{\max\left( \frac{4 \nu^2}{\theta} , \sqrt{  \frac{2 \nu^2 (F(\bbb{x}^0) - F(\bar{\bbb{x}}) - D)}{\theta } }\right)}{D}
\eeq
\noi After $c_{\bar{J}}$ iterations, Algorithm \ref{algo:main} becomes the proximal gradient method applied to the problem as in (\ref{eq:convex:subprob}). Therefore, the total number of iterations for finding a block-$k$ stationary point $N_1$ is bounded by:
\beq
N_1 &\overset{(a)}{\leq}& c_{\bar{J}}  + \frac{\max( \frac{4 \nu^2}{\theta} , \sqrt{  \frac{2 \nu^2 [F(\bbb{x}^{c_{\bar{J}}}) - F(\bar{\bbb{x}})]}{\theta } })}{\epsilon } \nn\\
&\overset{(b)}{\leq}& c_{\bar{J}}  + \frac{\max( \frac{4 \nu^2}{\theta} , \sqrt{  \frac{2 \nu^2 [ F(\bbb{x}^0) - F(\bar{\bbb{x}}) - D ]}{\theta } })}{\epsilon } \nn\\
&\overset{}{=}& \left( \frac{\bar{J}}{D} + \frac{1}{\epsilon}\right) \times  \max\left( \frac{4 \nu^2}{\theta} , \sqrt{  \frac{2 \nu^2 (F(\bbb{x}^0) - F(\bar{\bbb{x}}) - D)}{\theta } }\right)   \nn
\eeq
\noi where step $(a)$ uses the fact that the total number of iterations for finding a stationary point after $\bbb{x}_{c_{\bar{J}}}$ is upper bounded by ${\max( \frac{4 \nu^2}{\theta} , \sqrt{  \frac{2 \nu^2 [F(\bbb{x}^{c_{\bar{J}}}) - F(\bar{\bbb{x}})]}{\theta } })}/{\epsilon }$; step $(b)$ uses (\ref{eq:J:bound00}) and $j\geq 1$.

(ii) We now discuss the case when $f(\cdot)$ is strongly convex. Using part (ii) in Proposition \ref{the:sparse:bound:strong:yn}, we have:
\beq\label{eq:J:bound2}
c_j-c_{j-1} \leq \log_{\frac{\varpi}{1+\varpi}} \frac{D}{F(\bbb{x}^0)-F(\bar{\bbb{x}})} ,
\eeq
\noi Summing up the inequality (\ref{eq:J:bound2}) for $j = 1,2,...,\bar{J}$, we obtain the following results:
\beq \label{eq:jjjjjjjjjj}
c_{{\bar{J}}} \leq \log_{\frac{\varpi}{1+\varpi}} \left( \frac{D^{\bar{J}} }{(F(\bbb{x}^0)-F(\bar{\bbb{x}}))^{\bar{J}}} \right) = \bar{J} \log_{\frac{\varpi}{1+\varpi}} \left( \frac{D}{(F(\bbb{x}^0)-F(\bar{\bbb{x}}))} \right) \nn
\eeq
\noi Therefore, the total number of iterations $N_2$ is bounded by:
\beq
N_{2} &\overset{(a)}{\leq}& c_{\bar{J}}  + \log_{\frac{\varpi}{1+\varpi}} \left(\frac{\epsilon}{F(\bbb{x}^{c_{\bar{J}}}) - F(\bar{\bbb{x}})}\right) \nn\\
 &\overset{(b)}{\leq}& c_{\bar{J}}  + \log_{\frac{\varpi}{1+\varpi}} \left(\frac{\epsilon}{F(\bbb{x}^{0}) - D - F(\bar{\bbb{x}})}\right) \nn
\eeq
\noi where step $(a)$ uses the fact that the total number of iterations for finding a stationary point after $\bbb{x}_{c_{\bar{J}}}$ is upper bounded by $\log_{\frac{\varpi}{1+\varpi}} \left(\frac{\epsilon}{F(\bbb{x}^{c_{\bar{J}}}) - F(\bar{\bbb{x}})}\right)$; step $(b)$ uses (\ref{eq:J:bound00}) that $0\leq  F(\bbb{x}^0) - F(\bar{\bbb{x}}) -  t \times D$ and $t\geq 1$.

\end{proof}
\end{theorem}

\section{Convergence Rate for Sparsity Constrained Optimization} \label{sect:rate:sparse:c}

We now prove the convergence rate of Algorithm \ref{algo:main} for sparse constrained optimization with $h \triangleq f_{\text{sparse-c}}$. We results are based on the strongly convex and Lipschitz continuity of the objective function. We naturally derive the following theorem.

\setcounter{theorem}{3}
\begin{theorem}
\textbf{Proof of Convergence Rate when $h \triangleq f_{\text{sparse-c}}$.} Let $f(\cdot)$ be an $s$-strongly convex function. We assume that $f(\cdot)$ is Lipschitz continuous such that $\forall t,~\|\nabla f(\bbb{x}^t)\|_{2}^2\leq M$ for some constant $M$. We have the following results:
\beq \label{eq:sparse:cons}
\E[F(\bbb{x}^{t}) - F(\bar{\bbb{x}})] &\leq&   (F(\bbb{x}^{0}) - F(\bar{\bbb{x}}))\alpha^t +  \frac{M}{2\theta} \frac{\alpha}{1-\alpha},~\nn\\
\E[\frac{s}{4}\|\bbb{x}^{t+1} - \bar{\bbb{x}}\|^2_2] &\leq&   \frac{n  2\theta}{k}      (F(\bbb{x}^{0}) - F(\bar{\bbb{x}}))\alpha^t + \frac{n}{k}\frac{M}{1-\alpha},~\text{with}~\alpha \triangleq \frac{\frac{n\theta}{k s}}{1+\frac{n\theta}{ks}}.\nn
\eeq


\begin{proof}

(i) First of all, we define the zero set and nonzero set as follows:
\beq
S\triangleq \{i~|~i\in B,~\bbb{x}_{i}^{t+1} \neq 0\},~Z\triangleq \{i~|~i\in B,~\bbb{x}_{i}^{t+1} = 0\},~\nn
\eeq
\noi Using the optimality of $\bbb{x}^{t+1}$ for the subproblem, we obtain
\beq \label{eq:optimality:xk1}
(\nabla f(\bbb{x}^{t+1}))_{S} + \theta (\bbb{x}^{t+1}_S-\bbb{x}^{t}_S) = 0 ~\Rightarrow~(\nabla f(\bbb{x}^{t+1}))_{S}/\theta +  \bbb{x}^{t+1}_S  =  \bbb{x}^{t}_S
\eeq

\noi We derive the following inequalities:
\beq \label{eq:strongly:convex:sparse:rate}
&&  \E[f(\bbb{x}^{t+1}) - f(\bar{\bbb{x}})] \nn\\
&\overset{(a)}{\leq}&\E[ \la \bbb{x}^{t+1} - \bar{\bbb{x}},~\nabla f(\bbb{x}^{t+1}) \ra - \frac{s}{2}\|\bbb{x}^{t+1} - \bar{\bbb{x}}\|_2^2 ]\nn\\
&\overset{(b)}{=} & \frac{n}{k} \cdot \E[ \la \bbb{x}_{B}^{t+1} - \bbb{x}_{B}^*,~(\nabla f(\bbb{x}^{t+1}))_{B}\ra - \frac{s}{2}\|\bbb{x}^{t+1}_B - \bar{\bbb{x}}_B\|_{2}^2 ]   \nn\\
&\overset{(c)}{\leq} & \frac{n}{k} \cdot \frac{s}{2} \cdot \E[\|(\nabla f(\bbb{x}^{t+1}))_{B}/s\|_{2}^2 ]  \nn\\
&\overset{(d)}{=} & \frac{n}{k} \cdot \frac{s}{2} \cdot \left( \E[\|(\nabla f(\bbb{x}^{t+1}))_{S}/s\|_{2}^2 ] + \E[\|(\nabla f(\bbb{x}^{t+1}))_{Z}/s\|_{2}^2 ]\right)  \nn\\
& \overset{(e)}{\leq} & \E[ \frac{n}{k}\frac{s}{2} \cdot [ \| \theta(\bbb{x}_S^t - \bbb{x}^{t+1}_S)/s\|_2^2  +\|(\nabla f(\bbb{x}^{t+1}))_Z/s\|_2^2 ] ]\nn\\
& \overset{(f)}{\leq} &   \E[ \frac{n}{k}\frac{\theta}{s}\cdot \frac{\theta}{2}\| \bbb{x}^t - \bbb{x}^{t+1}\|_2^2  + \frac{n}{k}\frac{s}{2}\| (\nabla f(\bbb{x}^{t+1}))_Z/s\|_2^2 ]\nn\\
& \overset{(g)}{\leq} & \E[   \frac{n}{k}\frac{\theta}{s}\cdot [f(\bbb{x}^t)-f(\bbb{x}^{t+1})]  + \frac{n}{k}\frac{1}{2s}\|(\nabla f(\bbb{x}^{t+1}))_Z\|_2^2 ]\nn\\
& \overset{(h)}{\leq} & \E[ \frac{n}{k}\frac{\theta}{s}\cdot(f(\bbb{x}^t)-f(\bar{\bbb{x}}))-\frac{n}{k}\frac{\theta}{s}\cdot(f(\bbb{x}^{t+1})-f(\bar{\bbb{x}}))  + \frac{n}{k}\frac{1}{2s}\cdot M]
\eeq
\noi where step $(a)$ uses the strongly convexity of $f(\cdot)$; step $(b)$ uses the fact that the working set $B$ is selected with $\frac{k}{n}$ probability; step $(c)$ uses the inequality that $\la \bbb{x},~\bbb{a}\ra - \frac{s}{2}\|\bbb{x}\|_2^2 = \frac{s}{2}\|\bbb{a}/s\|_2^2-\frac{s}{2}\|\bbb{x}-\bbb{a}/s\|_2^2\leq\frac{s}{2}\|\bbb{a}/s\|_2^2$ for all $\bbb{a},~\bbb{x}$; step $(d)$ uses the fact that $B = S\cup Z$, step $(e)$ uses (\ref{eq:optimality:xk1}); step $(f)$ uses $\|\bbb{x}_S\|_2^2\leq \|\bbb{x}\|_2^2$ for all $\bbb{x}$; step $(g)$ uses the sufficient decrease condition that $\frac{\theta}{2}\|\bbb{x}^{t+1}-\bbb{x}^{t}\|^2 \leq F(\bbb{x}^t) - F(\bbb{x}^{t+1})$; step $(h)$ uses the Lipschitz continuity of $f(\cdot)$ that $\|\nabla f(\bbb{x}^{t+1})\|_2^2\leq M,~\forall t$.

\noi From (\ref{eq:strongly:convex:sparse:rate}), we have the following inequalities:
\beq
&&\E[(1+\frac{n\theta}{ks}) (f(\bbb{x}^{t+1}) - f(\bar{\bbb{x}}))] \leq  \E[\frac{n \theta}{ks}\cdot(f(\bbb{x}^t)-f(\bar{\bbb{x}})) + \frac{n M}{2ks} ] \nn\\
&&\E[f(\bbb{x}^{t+1}) - f(\bar{\bbb{x}})] \leq \E[ \alpha (f(\bbb{x}^t)-f(\bar{\bbb{x}})) + \frac{\frac{n}{2ks}}{\frac{n}{k}\frac{\theta}{s}} \alpha M ]\nn\\
&& \E[ f(\bbb{x}^{t+1}) - f(\bar{\bbb{x}})] \leq \E[ \alpha (f(\bbb{x}^t)-f(\bar{\bbb{x}})) +  \frac{\alpha M}{2\theta}] \nn
\eeq

\noi Solving this recursive formulation, we have:
\beq \label{eq:convergence:rate:sparse:conc:1}
\E[ f(\bbb{x}^{t}) - f(\bar{\bbb{x}}) ] &\leq&  \E[ \alpha^t (f(\bbb{x}^{0}) - f(\bar{\bbb{x}}))] +  M \sum_{i=1}^t \alpha^i \nn\\
&=& \E[\alpha^t (f(\bbb{x}^{0}) - f(\bar{\bbb{x}}))] +  \frac{M \alpha}{2\theta} \cdot \frac{1-\alpha^t}{1-\alpha}\nn\\
&\leq & \E[\alpha^t (f(\bbb{x}^{0}) - f(\bar{\bbb{x}}))] +  \frac{M }{2\theta} \cdot \frac{\alpha}{1-\alpha}\nn
\eeq
\noi Since $\bbb{x}^t$ is always a feasible solution for all $t=1,2,...\infty$, we have $F(\bbb{x}^t)=f(\bbb{x}^t)$. Therefore, we obtain (\ref{eq:sparse:cons}).

\noi (ii) We now prove the second part of this theorem. First, we derive the following inequalities:
\beq \label{eq:bound:xt1xt}
\E[\|\bbb{x}^{t+1}-\bbb{x}^{t}\|_2^2] &\overset{(a)}{\leq} & \E[\frac{2}{\theta} \left(F(\bbb{x}^t) - F(\bbb{x}^{t+1}) \right) ]\nn\\
&\overset{(b)}{\leq}& \E[\frac{2}{\theta} \left( F(\bbb{x}^t) - F(\bar{\bbb{x}}) \right)] \nn\\
&\overset{(c)}{\leq}& \E[ \frac{2}{\theta} \alpha^t (f(\bbb{x}^{0}) - f(\bar{\bbb{x}}))] +  \frac{M }{\theta^2} \cdot \frac{\alpha}{1-\alpha}
\eeq
\noi where step $(a)$ uses the sufficient decrease condition in (\ref{eq:suff:dec}); step $(b)$ uses the fact that $F(\bar{\bbb{x}}) \leq F({\bbb{x}}^{t+1})$; step $(c)$ uses the result in (\ref{eq:convergence:rate:sparse:conc:1}).

Second, we have the following results:
\beq \label{eq:strongly:convex:sparse:rate:eq}
 \E[\frac{s}{2}\|\bbb{x}^{t+1} - \bar{\bbb{x}}\|_2^2] &\overset{(a)}{\leq}& \E[\la \bbb{x}^{t+1} - \bar{\bbb{x}},~\nabla f(\bbb{x}^{t+1}) \ra  + f(\bar{\bbb{x}})  - f(\bbb{x}^{t+1})]  \nn\\
&\overset{(b)}{\leq}& \E[ \la \bbb{x}^{t+1} - \bar{\bbb{x}},~\nabla f(\bbb{x}^{t+1}) \ra ] \nn\\
&\overset{(c)}{\leq}&  \E[\|\bbb{x}^{t+1} - \bar{\bbb{x}}\| \cdot \|\nabla f(\bbb{x}^{t+1})\|]
\eeq
\noi where step $(a)$ uses the strongly convexity of $f(\cdot)$; step $(b)$ uses the fact that $f(\bar{\bbb{x}}) \leq f({\bbb{x}}^{t+1})$; step $(c)$ uses the Cauchy-Schwarz inequality.

From (\ref{eq:strongly:convex:sparse:rate:eq}), we further have the following results:
\beq \label{eq:strongly:convex:sparse:eq2}
\E[\frac{s}{4}\|\bbb{x}^{t+1} - \bar{\bbb{x}}\|^2_2] &\overset{(a)}{\leq}&\E[ \|\nabla f(\bbb{x}^{t+1})\|_2^2  ]\nn\\
&\overset{}{ = } & \frac{n}{k} \E[  \|\nabla_B f(\bbb{x}^{t+1})\|_2^2  ]\nn\\
&\overset{(b)}{ = } &  \frac{n}{k}  \E[  \|\nabla_S f(\bbb{x}^{t+1})\|_2^2 +  \|\nabla_Z f(\bbb{x}^{t+1})\|_2^2   ]\nn\\
&\overset{(c)}{ \leq } &   \frac{n}{k} \E[  \theta^2\| \bbb{x}^{t+1}_S  - \bbb{x}^{t}_S \|_2^2] +   \frac{n}{k} M   \nn\\
&\overset{(d)}{ = } &   \frac{n}{k}  \E[ {2 \theta} \alpha^t (f(\bbb{x}^{0}) - f(\bar{\bbb{x}})) +   \frac{M}{1-\alpha} ]  \nn
\eeq
\noi where step $(a)$ uses the strongly convexity of $f(\cdot)$; step $(b)$ uses the fact that $B = S \cup Z$; step $(c)$ uses the assumption that $\|\nabla f(\bbb{x}^t)\|_2^2\leq M$ for all $\bbb{x}^t$ and the optimality of $\bbb{x}^{t+1}$ in (\ref{eq:optimality:xk1}); step $(d)$ uses (\ref{eq:bound:xt1xt}). Therefore, we finish the proof of this theorem.

\end{proof}

\end{theorem}

\clearpage
\bibliography{my}

\end{document}